\documentclass[10pt,leqno,draft,twoside]{article}

\usepackage{amsmath,amsthm,amssymb}
\numberwithin{equation}{section}

\makeatletter
\newif\ify@autoscale \y@autoscaletrue \def\Yautoscale#1{\ifnum #1=0
  \y@autoscalefalse\else\y@autoscaletrue\fi}
\newdimen\y@b@xdim
\newdimen\y@boxdim \y@boxdim=13pt
\def\y@style{\displaystyle} 
\def\Yboxdim#1{\y@autoscalefalse\y@boxdim=#1}
\newdimen\y@linethick    \y@linethick=.3pt
\def\Ylinethick#1{\y@linethick=#1}
\newskip\y@interspace \y@interspace=0ex plus 0.3ex
\def\Yinterspace#1{\y@interspace=#1}
\newif\ify@vcenter   \y@vcentertrue	
\def\Yvcentermath#1{\ifnum #1=0 \y@vcenterfalse\else\y@vcentertrue\fi}
\newif\ify@stdtext   \y@stdtextfalse
\def\Ystdtext#1{\ifnum #1=0 \y@stdtextfalse\else\y@stdtexttrue\fi}
\def\y@vr{\vrule height0.8\y@b@xdim width\y@linethick depth 0.2\y@b@xdim}
\def\y@emptybox{\y@vr\hbox to \y@b@xdim{\hfil}}
\def\y@abcbox#1{\y@vr\hbox to \y@b@xdim{\hfil#1\hfil}}
\def\y@mathabcbox#1{\y@vr\hbox to \y@b@xdim{\hfil${\y@style#1}$\hfil}} 
\def\y@setdim{%
  \ify@autoscale%
   \setbox1=\hbox{A}\y@b@xdim=1.6\ht1 \setbox1=\hbox{}\box1%
  \else\y@b@xdim=\y@boxdim \advance\y@b@xdim by -2\y@linethick
  \fi}
\newcount\y@counter
\newif\ify@islastarg
\def\y@lastargtest#1,#2 {\if\space #2 \y@islastargtrue
  \else\y@islastargfalse\fi}
\def\y@emptyboxes#1{\y@counter=#1\loop\ifnum\y@counter>0
  \advance\y@counter by -1 \y@emptybox\repeat}
\def\y@nelineemptyboxes#1{%
  \vbox{%
    \hrule height\y@linethick%
    \hbox{\y@emptyboxes{#1}\y@vr}
    \hrule height\y@linethick}\vspace{-\y@linethick}}
\def\yng(#1){%
  \y@setdim%
  \hspace{\y@interspace}%
  \ifmmode\ify@vcenter\vcenter\fi\fi{%
  \y@lastargtest#1,
  \vbox{\offinterlineskip
    \ify@islastarg
     \y@nelineemptyboxes{#1}
    \else
     \y@ungempty(#1)
    \fi}}\hspace{\y@interspace}}
\def\y@ungempty(#1,#2){%
  \y@nelineemptyboxes{#1}
  \y@lastargtest#2,
  \ify@islastarg
   \y@nelineemptyboxes{#2}
  \else
   \y@ungempty(#2)
  \fi}
\def\y@nelettertest#1#2. {\if\space #2 \y@islastargtrue
  \else\y@islastargfalse\fi}
\def\y@abcboxes#1#2.{%
  \ify@stdtext\y@abcbox#1\else\y@mathabcbox#1\fi%
  \y@nelettertest #2.
  \ify@islastarg\unskip%
   \ify@stdtext\y@abcbox{#2}\else\y@mathabcbox{#2}\fi%
  \else\y@abcboxes#2.\fi}
\def\y@m@veright@ifskew#1{}
\def\y@nelineabcboxes#1{%
  \y@nelettertest #1.
  \ify@islastarg
   \y@m@veright@ifskew{#1}
    \vbox{
      \hrule height\y@linethick%
      \hbox{\ify@stdtext\y@abcbox#1\else\y@mathabcbox#1\fi\y@vr}
      \hrule height\y@linethick}\vspace{-\y@linethick}
  \else
   \y@m@veright@ifskew{#1}
    \vbox{
      \hrule height\y@linethick%
      \hbox{\y@abcboxes #1.\y@vr}%
      \hrule height\y@linethick}\vspace{-\y@linethick}
  \fi}
\def\young(#1){%
  \y@setdim%
  \hspace{\y@interspace}%
  \y@lastargtest#1,
  \ifmmode\ify@vcenter\vcenter\fi\fi{%
  \vbox{\offinterlineskip
    \ify@islastarg\y@nelineabcboxes{#1}%
    \else\y@ungabc(#1)%
    \fi}}\hspace{\y@interspace}}
\def\y@ungabc(#1,#2){%
  \y@nelineabcboxes{#1}%
  \y@lastargtest#2,
  \ify@islastarg\y@nelineabcboxes{#2}%
  \else\y@ungabc(#2)%
  \fi}

\def\minitab(#1){%
\y@autoscalefalse%
\def\y@style{\scriptscriptstyle}%
\y@boxdim=6pt%
\young(#1)%
\y@autoscaletrue%
\def\y@style{\displaystyle}%
\y@boxdim=13pt}
\def\tinytab(#1){%
\y@autoscalefalse%
\def\y@style{\scriptscriptstyle}%
\y@boxdim=5pt%
\young(#1)%
\y@autoscaletrue%
\def\y@style{\displaystyle}%
\y@boxdim=13pt}

\def\miniyng(#1){%
\y@autoscalefalse%
\y@boxdim=6pt%
\yng(#1)%
\y@autoscaletrue%
\y@boxdim=13pt}
\def\tinyyng(#1){%
\y@autoscalefalse%
\y@boxdim=5pt%
\yng(#1)%
\y@autoscaletrue%
\y@boxdim=13pt}

\makeatother

\newtheorem{thm}{Theorem}[section]
\newtheorem{prop}[thm]{Proposition}
\newtheorem{lem}[thm]{Lemma}
\newtheorem{cor}[thm]{Corollary}

\newtheorem*{conjA}{Conjecture A}
\newtheorem*{conjAp}{Conjecture $\text{A\!}'$}
\newtheorem*{conjB}{Conjecture B}
\newtheorem*{conjAB}{Conjecture AB}

\theoremstyle{definition}

\newtheorem{example}[thm]{Example}

\newtheorem*{ackn}{Acknowledgement}

\theoremstyle{remark}
\newtheorem{remark}[thm]{Remark}

\newcommand{\eoe}{\hfill{\rule{0.7em}{0.7em}}}
\newenvironment{ex}{\begin{example}}
{\end{example}}
\newenvironment{rem}{\begin{remark}}
{\eoe\end{remark}}

\newenvironment{keywords}{\smallskip\noindent{\bfseries Keywords:}}{}
\newenvironment{MSC}
{\smallskip\noindent{\bfseries 2000 Mathematical Subject Classification:}}{}
\newenvironment{dedication}{\begin{center}\itshape}{\end{center}}

\newcommand{\card}[1]{\left|#1\right|}	
\newcommand{\abs}[1]{\left|#1\right|}	
\newcommand{\sym}[1]{\mathfrak{S}_{#1}}	
\newcommand{\pp}[1]{{\upshape(}#1{\upshape)}}	

\newcommand{\deq}{:=}	
\newcommand{\inv}[1]{\ell(#1)}	
\newcommand{\cyc}[2][n]{\nu_{#1}(#2)}
\newcommand{\tca}{\curvearrowleft}	

\newcommand{\blf}[2]{\langle#1,\,#2\rangle}	
\newcommand{\q}[3]{q_{#1}^{#2}(#3)}
\newcommand{\cop}{\Delta}	
\newcommand{\e}{\varepsilon}	

\newcommand{\Z}{\mathbb{Z}}	
\newcommand{\C}{\mathbb{C}} 

\newcommand{\ve}{\boldsymbol{e}}	
\newcommand{\vi}{\boldsymbol{i}}
\newcommand{\vj}{\boldsymbol{j}}
\newcommand{\vk}{\boldsymbol{k}}

\newcommand{\uF}{\widetilde{F}}	
\newcommand{\uC}{\widetilde{C}}	
\newcommand{\uQ}{\widetilde{Q}}

\let \det = \relax
\DeclareMathOperator{\det}{det}	
\DeclareMathOperator{\Mat}{Mat}	
\DeclareMathOperator{\End}{End}	
\DeclareMathOperator{\rank}{rank}	
\DeclareMathOperator{\diag}{diag}	
\DeclareMathOperator{\STab}{STab}	
\DeclareMathOperator{\per}{per}	
\DeclareMathOperator{\Sing}{Sing}	
\DeclareMathOperator{\prSing}{prSing}
\DeclareMathOperator{\wt}{wt}	

\newcommand{\singular}[2][q]{\Sing_{#2,#1}}
\newcommand{\psingular}[2][q]{\prSing_{#2,#1}}
\newcommand{\subsingular}[3][q]{\Sing^{#3}_{#1}}
\newcommand{\singularwt}[3][q]{\Sing_{#2,#1}\!\left(#3\right)}
\newcommand{\persingular}[2][q]{\Sing^{\mathrm{per}}_{#2,#1}}

\newcommand{\liegl}{\mathfrak{gl}}	
\newcommand{\A}[1][n]{\mathcal{A}(\Mat_{#1})}	
\newcommand{\U}[1][n]{\mathcal{U}(\liegl_{#1})}	
\newcommand{\Aq}[1][n]{\mathcal{A}_q(\Mat_{#1})}	
\newcommand{\Uq}[1][n]{\mathcal{U}_q(\liegl_{#1})}	
\newcommand{\Hq}[1][n]{\mathcal{H}_q(\sym{#1})}	

\newcommand{\Umod}[3][q]{\boldsymbol{E}_{\!#2,#1}^{#3}}	
\newcommand{\Hmod}[3][q]{\boldsymbol{M}_{\!#2,#1}^{#3}}	

\newcommand{\linspan}[3][q]{\boldsymbol{L}_{#2,#1}^{(#3)}}
\newcommand{\cycmod}[3][q]{\boldsymbol{V}_{\!\!#2,#1}^{(#3)}}
\newcommand{\cycmodp}[3][q]{\boldsymbol{U}_{\!#2,#1}^{(#3)}}
\newcommand{\UHcycmod}[3][q]{\mathcal{V}_{#2,#1}^{(#3)}}

\newcommand{\lattice}[1][n]{\mathcal{L}_{#1}}	
\newcommand{\ulattice}[1][n]{\widetilde{\mathcal{L}}_{#1}}	
\newcommand{\domlattice}[1][n]{{\ulattice[#1]}^{\mathrm{dom}}}	

\newcommand{\hwv}[2][\alpha]{v^{(#1)}(#2)}	

\newcommand{\yngsym}[1][q]{\mathbb{E}_{#1}}	

\newcommand{\mult}[3][q]{m_{#1}^{#2}\!\left(#3\right)}	
\newcommand{\multp}[3][q]{m_{#1}^{#2}\!\left(#3\right)_{\mathrm{per}}}

\newcommand{\Det}[1][\alpha]{\det_q^{(#1)}}	
\newcommand{\Per}[1][\alpha]{\per_q^{(#1)}}	
\newcommand{\qad}[2][\alpha]{D_q^{(#1)}(#2)}
\newcommand{\qap}[2][\alpha]{P_q^{(#1)}(#2)}

\newcommand{\rhoCn}{\rho_{\lower0.5ex\hbox{$\scriptstyle\C^n$}}}
\newcommand{\rhoCnn}{\rhoCn^{\otimes n}}
\newcommand{\rhogl}[1][n]{\rho_{\!\lower0.5ex\hbox{$\scriptstyle\liegl_{#1}$}}\!}
\newcommand{\tX}{{}^{\raise0.35ex\hbox{$\scriptstyle t$}}\kern-0.2em X}

\newcommand{\LRarrow}[2]{%
\ \underset{#2}{\overset{#1}%
{\ooalign{\raise0.45ex\hbox{$\longrightarrow$}%
\crcr\raise-0.45ex\hbox{$\longleftarrow$}}}}\ }
\newcommand{\Rarrow}[1]{\overset{#1}{\longrightarrow}}
\newcommand{\Larrow}[1]{\underset{#1}{\longleftarrow}}

\newcommand{\h}[1]
{\ifcase#1\or%
{1}\or{h_1}\or{h_2}\or{h_3}\or{h_1h_2h_1}\or{h_2h_3h_2}\or{h_1h_2h_3h_2h_1}\or%
{h_1h_2}\or{h_2h_1}\or{h_2h_3}\or{h_3h_2}\or{h_1h_2h_3h_2}\or{h_2h_3h_2h_1}\or%
{h_1h_2h_1h_3}\or{h_1h_3h_2h_1}\or{h_1h_2h_3}\or{h_3h_2h_1}\or{h_1h_3h_2}\or%
{h_2h_1h_3}\or{h_1h_2h_1h_3h_2}\or{h_2h_3h_2h_1h_2}\or{h_1h_3}\or%
{h_2h_1h_3h_2}\or{h_1h_2h_1h_3h_2h_1}\fi}

\newcommand{\set}[3][0]
{\ifcase#1%
\left\{#2\,;\,#3\right\}\or%
\bigl\{#2\,\big|\,#3\bigr\}\or%
\Bigl\{#2\,\Big|\,#3\Bigr\}\or%
\biggl\{#2\,\bigg|\,#3\biggr\}\or%
\Biggl\{#2\ \Bigg|\ #3\Biggr\}\fi}

\newcommand{\Dtheta}[1][q]{\vartheta_{#1}^{\mathrm{det}}}
\newcommand{\Ptheta}[1][q]{\vartheta_{#1}^{\mathrm{per}}}
\def\partitions{\mathcal{P}}
\def\piper{\pi_{\mathrm{per}}}

\makeatletter
\newcommand{\@newkakko}[4]
{\ifcase#1%
\left#3#2\right#4\or%
\bigl#3#2\bigr#4\or%
\Bigl#3#2\Bigr#4\or%
\biggl#3#2\biggr#4\or%
\Biggl#3#2\Biggr#4\else
#2\fi}

\newcommand{\kakko}[2][0]{\@newkakko{#1}{#2}{(}{)}}
\newcommand{\ckakko}[2][0]{\@newkakko{#1}{#2}{\{}{\}}}
\newcommand{\dkakko}[2][0]{\@newkakko{#1}{#2}{[}{]}}
\newcommand{\akakko}[2][0]{\@newkakko{#1}{#2}{\langle}{\rangle}}
\makeatother

\title{\bfseries Quantum $\alpha$-determinant cyclic modules of $\Uq$}
\author{Kazufumi Kimoto%
\thanks{Partially supported by Grant-in-Aid for Young Scientists (B) No.16740021.
(corresponding author)}\\
\small Department of Mathematical Science,
University of the Ryukyus,\\
\small Senbaru, Nishihara, Okinawa 903-0231, JAPAN
\\
\\
Masato Wakayama%
\thanks{Partially supported by Grant-in-Aid for Scientific Research (B) No. 15340012.}\\
\small Faculty of Mathematics, Kyushu University,
Hakozaki, Fukuoka 812-8518, JAPAN%
\footnote{Email addresses: \texttt{kimoto@math.u-ryukyu.ac.jp} (Kimoto),
\texttt{wakayama@math.kyushu-u.ac.jp} (Wakayama)}}

\begin{document}
\maketitle

\begin{dedication}
Dedicated to Takashi Ichinose on the occasion of his 65th birthday.
\end{dedication}

\begin{abstract}
As a particular one parameter deformation
of the quantum determinant,
we introduce
a quantum $\alpha$-determinant $\Det$
and study the $\Uq$-cyclic module generated by it:
We show that the multiplicity
of each irreducible representation
in this cyclic module is determined
by a certain polynomial called the $q$-content discriminant.
A part of the present result is
a quantum counterpart for the result of Matsumoto and Wakayama \cite{MW2005},
however, a new distinguished feature arises in our situation.
Specifically, we determine the degeneration
of the multiplicities for \emph{`classical'} singular points
and give a general conjecture for singular points
involving \emph{semi-classical} and \emph{quantum} singularities.
Moreover, we introduce a quantum $\alpha$-permanent $\Per$
and establish another conjecture
which describes a `reciprocity' between the multiplicities
of the irreducible summands of the cyclic modules
generated respectively by $\Det$ and $\Per$.


\begin{MSC}
20G42, 
20C08. 
\end{MSC}
\end{abstract}

\tableofcontents

\section{Introduction}

Let $\A$ be the associative $\C$-algebra
consisting of polynomial functions
on the set $\Mat_n$ of $n$ by $n$ matrices.
We denote by $x_{ij}$ the standard coordinate function on $\Mat_n$
with respect to the matrix unit $E_{ij}$.
The right translation of the general linear group $GL_n$ on $\A$
induces the representation $\rhogl$
of the enveloping algebra $\U$ of $\liegl_n=\liegl_n(\C)$ as
\begin{equation*}
\rhogl(e_{ij})=\sum_{k=1}^n x_{ki}\frac{\partial}{\partial x_{kj}}.
\end{equation*}
Here $\{e_{ij}\}_{1\le i,j\le n}$ is the standard basis
of $\liegl_n$ so that
$[e_{ij}, e_{kl}]=\delta_{jk}e_{il}-\delta_{li}e_{kj}$.
It is a very basic fact that the determinant
\begin{equation*}
\det X=\sum_{w\in\sym n}(-1)^{\inv w}x_{w(1)1}\dotsb x_{w(n)n}\in\A
\end{equation*}
of $X=\sum_{i,j}x_{ij}E_{ij}\in\Mat_n(\A)$
is an invariant of the action of $GL_n$.
In other words, we see that
$\rhogl(\U)\cdot\det X=\C\cdot\det X$.
We also have another distinguished element called the permanent
$\per X$ in $\A$ defined by
\begin{equation*}
\per X=\sum_{w\in\sym n}x_{w(1)1}\dotsb x_{w(n)n}.
\end{equation*}
Though $\per X$ is not an invariant of the action,
the cyclic module $\rhogl(\U)\cdot\per X$
gives an irreducible representation of $GL_n$
on the space of $n$-symmetric tensors of $\C^n$.


For a complex parameter $\alpha$,
the \emph{$\alpha$-determinant} is defined by
\begin{equation*}
\det^{(\alpha)}X=\sum_{w\in\sym n}\alpha^{n-\cyc w}
x_{w(1)1}x_{w(2)2}\dotsb x_{w(n)n}\in\A,
\end{equation*}
where $\cyc w$ is the number of cycles in $w\in\sym n$ \cite{ST2003}.
We notice that $\det^{(-1)}X=\det X$ and $\det^{(1)}X=\per X$.
Thus the $\alpha$-determinant interpolates
the determinant and permanent.
In the representation-theoretic point of view,
we can also understand that
the $\U$-cyclic module generated by
the $\alpha$-determinant interpolates the two irreducible representations;
the skew-symmetric tensor representation $\rhogl(\U)\cdot\det X$
and symmetric tensor representation $\rhogl(\U)\cdot\per X$.

Therefore it is natural to study the structure of
the interpolating module $V_n^{(\alpha)}\deq\rhogl(\U)\cdot\det^{(\alpha)}X$.
This is done in \cite{MW2005}.
The results are summarized as follows (see Section 4.2).
The module $V_n^{(\alpha)}$
is isomorphic to the tensor module
$(\C^n)^{\otimes n}\cong \bigoplus_{\lambda}E_\lambda^{\oplus f^\lambda}$
for all but finite exceptional values of $\alpha$.
Here $E_\lambda$ denotes the irreducible highest weight module of $\U$
of highest weight $\lambda$ and $f^\lambda$ the multiplicity of $E_\lambda$.
The isotypic component $E_\lambda^{\oplus f^\lambda}$ in $V_n^{(\alpha)}$
of highest weight $\lambda$ disappears
when $\alpha$ is a root of the (modified) \emph{content polynomial} $c_\lambda(x)$
\cite{Mac}.
Further, if we consider the cyclic module $\rhogl[2](\U[2])\cdot(\det^{(\alpha)}X)^k$
for positive integers $k$,
we see that
the disappearance of a subrepresentation
(from the cyclic module in a general position)
is described by
the \emph{Jacobi polynomial} with special parameters
determined by $k$ and the corresponding subrepresentation \cite{KMW2006(ip)}.
More precisely,
for each irreducible representation with highest weight $\lambda$ appearing
in $\rhogl[2](\U[2])\cdot(\det^{(\alpha)}X)^k$,
there is a polynomial $F_k^\lambda(x)$
such that the $\lambda$-isotypic component in $\rhogl[2](\U[2])\cdot(\det^{(\alpha)}X)^k$
is killed when $\alpha$ is a root of $F_k^\lambda(x)$,
and the polynomial $F_k^\lambda(x)$ is given as the Jacobi polynomial
whose parameters are explicitly determined by $k$ and $\lambda$.
Thus, we expect to find new families of polynomilas systematically
by considering the cyclic modules $\rhogl(\U)\cdot(\det^{(\alpha)}X)^k$ for $n\ge3$
as polynomials whose roots describe the degeneration of
the module $\rhogl(\U)\cdot(\det^{(\alpha)}X)^k$.

The point of the story is also that
only one element $\det^{(\alpha)}X$ generates
various irreducible representations
with emphasizing that special polynomials describe
the degeneration of the module $\rhogl(\U)\cdot\det^{(\alpha)}X$.
This is quite a \emph{contrast} to the standard representation theory
in which we construct the irreducible $\U$-modules
by utilizing various minor determinants as highest weight vectors and
we get various special functions as matrix coefficients of them.

The situation allows us to propose a strategy
for discovering new special polynomials
as polynomials which controls the structure of cyclic $\U$-modules
generated by the $\alpha$-determinants.
Since there are several special functions such as the
Jacobi (big/little) $q$-polynomials which we obtain as matrix coefficients
of irreducible representations of \emph{quantum groups},
it is natural to quantize the situation described above to get wider class of
special polynomials.




We take the natural representation $\rho$
of the quantum enveloping algebra $\Uq$
on the quantum matrix algebra $\Aq$.
We introduce a natural quantization $\Det X\in\Aq$
of the $\alpha$-determinant $\det^{(\alpha)}$
as a particular deformation
of the quantum determinant $\det_q X$ in \cite{J1986} by the formula
\begin{equation*}
\Det X \deq \sum_{w\in\sym{n}}
\alpha^{n-\cyc w}q^{\inv w}x_{w(1)1}\dotsb x_{w(n)n}
\in\Aq.
\end{equation*}
We call this element $\Det$ \emph{quantum $\alpha$-determinant}.
We notice that the quantum $\alpha$-determinant
$\Det[-1]$ for $\alpha=-1$
is nothing but
the quantum determinant $\det_q$.

In the present paper,
as a beginning of the study,
we treat the irreducible decomposition
of the cyclic module
$\cycmod n\alpha\deq\rho(\Uq)\cdot\Det$
and show that
there are finite number of values called \emph{singular values}
such that
the structure of $\cycmod n\alpha$ changes drastically
if $\alpha$ is one of such values,
and the singular values are actually described as roots of some polynomials
called \emph{$q$-content discriminants}.

We now briefly sketch the contents of the paper below.

The basic conventions on the quantum matrix algebra $\Aq$,
the quantum enveloping algebra $\Uq$ (as a quantum group)
and the Iwahori-Hecke algebra $\Hq$
are collected briefly in Section \ref{Preliminaries}.

In Section \ref{QAD},
we define the quantum $\alpha$-determinant $\Det$
and investigate several basic properties
of the module $\cycmod n\alpha$.
As a starting point,
we show in Proposition \ref{prop:ird_generic} that
$\cycmod n\alpha$ is equivalent
to the tensor product module $(\C^n)^{\otimes n}$
for all but finite $\alpha$.
Namely, for almost all values $\alpha$,
the multiplicity $\mult\lambda\alpha$
of the highest weight module $\Umod n{\lambda}$
corresponding to a highest weight ($=$ partition) $\lambda$
in $\cycmod n\alpha$
is equal to the number $f^\lambda$
of standard tableaux with shape $\lambda$.
If $\alpha$ is one of those finite exceptions,
that is, if $\mult\lambda\alpha<f^\lambda$ holds,
we call it a singular point (or value) as we mentioned above.
To describe highest weight vectors
of the irreducible factors of the decomposition of $\cycmod n\alpha$
in terms of the quantum $\alpha$-determinants,
we employ the $q$-Young symmetrizers studied in \cite{G1986}.

The degeneration of the cyclic module $\cycmod n\alpha$
is discussed in Section \ref{Main}.
For each $\lambda$,
a $f^\lambda\times f^\lambda$ matrix called a $q$-content transition matrix
and a polynomial in $\alpha$ called a $q$-content discriminant
are introduced.
The zeros of the $q$-content discriminants are
the singular values.
In contrast with the classical theory developed in \cite{MW2005},
the $q$-Young symmetrizers can not give us enough information
about the zeros of $q$-content discriminants,
and the explicit description of the zeros of these polynomials
seems to be a far reaching problem.
This is the most difficult point in the present study
which we have never encountered
in the classical situation \cite{MW2005}.

In the classical theory,
the content transition matrix is a scalar one
whose scalar is given by the so-called content polynomial (see \cite{Mac}).
It follows hence that
only $\pm\frac1k\;(1\leq k< n)$ are the singular points
and the degeneration $m_1^\lambda(\pm\frac1k)<f^\lambda$
implies the vanishing $m_1^\lambda(\pm\frac1k)=0$.
In the quantization,
the singular points $-1,-\frac12,\dots,-\frac1{n-1}$
remains singular,
whereas the points $1,\frac12,\ldots,\frac1{n-1}$ themselves
are no longer singular but are $q$-deformed,
and even new values (depending on $q$) other than
$2n-2$ singular values above come up as extra singular points.
We call the first member of singular points \emph{classical},
the second one \emph{semi-classical}
and the last one \emph{quantum} respectively.
We devote ourselves to investigate such singular points
in the latter half of the section
and give one conjecture
concerning the multiplicity degeneration for singular values
(Conjecture A, see also Theorem \ref{thm:classical sing}).
Indeed, for our quantum case,
it occurs that $0<\mult\lambda\alpha<f^\lambda$
for a quantum singular point $\alpha$.
We furthermore translate this conjecture
into the framework of the bimodule of $\Uq$ and $\Hq$
in terms of Schur-Weyl duality \cite{G1986, J1986}.
We also provide explicit calculations
of some $q$-content discriminants for readers' help.
It would be also interesting to study the relation
between the multiplicity $\mult\lambda\alpha$
and the multiplicity of $\alpha$ as a root of the corresponding
$q$-content discriminant.
We treat this subject in the future study.

In the last section,
we introduce a notion of quantum $\alpha$-permanent $\Per=\det_{q^{-1}}^{(-\alpha)}$
and discuss shortly the $\Uq$-cyclic module $\cycmodp n\alpha$
generated by $\Per$ through examples.
In the classical case,
as we mentioned above,
the $\U$-cyclic module generated by
the $\alpha$-determinant interpolates the two irreducible representations
$\rhogl(\U)\cdot\det^{(-1)}X$ and $\rhogl(\U)\cdot\det^{(1)} X$.
Whereas, in the quantum situation,
the cyclic module $\cycmod n\alpha$ is irreducible
(skew-symmetric tensor representation)
at $\alpha=-1$ but not irreducible at $\alpha=1$,
as well as
the cyclic module $\cycmodp n\alpha$ is irreducible
(symmetric tensor representation)
at $\alpha=-1$ but not irreducible at $\alpha=1$.
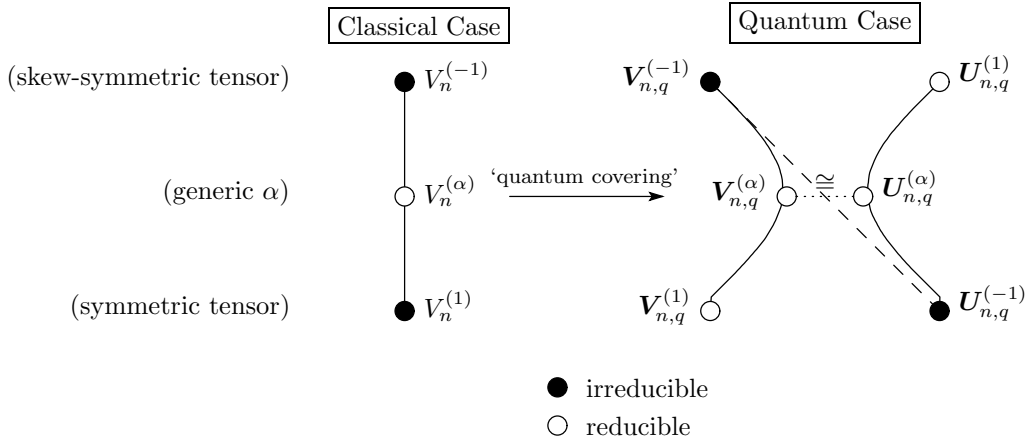
\begin{flushleft}
\medskip
\unitlength 0.1in
\begin{picture}( 55.7000, 22.2000)( -2.7000,-24.5000)
\special{pn 8}%
\special{pa 2400 1800}%
\special{pa 2400 600}%
\special{fp}%
\special{pn 8}%
\special{sh 1.000}%
\special{ar 2400 600 50 50  0.0000000 6.2831853}%
\special{pn 8}%
\special{sh 1.000}%
\special{ar 2400 1800 50 50  0.0000000 6.2831853}%
\special{pn 8}%
\special{sh 0}%
\special{ar 2400 1200 50 50  0.0000000 6.2831853}%
\special{pn 8}%
\special{pa 4000 600}%
\special{pa 4028 628}%
\special{pa 4054 654}%
\special{pa 4080 680}%
\special{pa 4106 708}%
\special{pa 4132 734}%
\special{pa 4156 760}%
\special{pa 4180 786}%
\special{pa 4204 814}%
\special{pa 4228 840}%
\special{pa 4250 866}%
\special{pa 4270 894}%
\special{pa 4290 920}%
\special{pa 4308 946}%
\special{pa 4326 974}%
\special{pa 4340 1000}%
\special{pa 4354 1026}%
\special{pa 4368 1054}%
\special{pa 4378 1080}%
\special{pa 4386 1106}%
\special{pa 4394 1134}%
\special{pa 4398 1160}%
\special{pa 4400 1186}%
\special{pa 4400 1212}%
\special{pa 4398 1240}%
\special{pa 4394 1266}%
\special{pa 4388 1292}%
\special{pa 4378 1320}%
\special{pa 4368 1346}%
\special{pa 4356 1372}%
\special{pa 4342 1400}%
\special{pa 4326 1426}%
\special{pa 4310 1452}%
\special{pa 4292 1480}%
\special{pa 4272 1506}%
\special{pa 4250 1532}%
\special{pa 4228 1560}%
\special{pa 4206 1586}%
\special{pa 4182 1612}%
\special{pa 4158 1638}%
\special{pa 4134 1666}%
\special{pa 4108 1692}%
\special{pa 4082 1718}%
\special{pa 4056 1746}%
\special{pa 4028 1772}%
\special{pa 4002 1798}%
\special{pa 4000 1800}%
\special{sp}%
\special{pn 8}%
\special{pa 5200 600}%
\special{pa 5174 628}%
\special{pa 5148 654}%
\special{pa 5122 680}%
\special{pa 5096 708}%
\special{pa 5070 734}%
\special{pa 5044 760}%
\special{pa 5020 786}%
\special{pa 4996 814}%
\special{pa 4974 840}%
\special{pa 4952 866}%
\special{pa 4930 894}%
\special{pa 4912 920}%
\special{pa 4892 946}%
\special{pa 4876 974}%
\special{pa 4860 1000}%
\special{pa 4846 1026}%
\special{pa 4834 1054}%
\special{pa 4824 1080}%
\special{pa 4814 1106}%
\special{pa 4808 1134}%
\special{pa 4804 1160}%
\special{pa 4800 1186}%
\special{pa 4800 1212}%
\special{pa 4802 1240}%
\special{pa 4808 1266}%
\special{pa 4814 1292}%
\special{pa 4822 1320}%
\special{pa 4832 1346}%
\special{pa 4846 1372}%
\special{pa 4860 1400}%
\special{pa 4874 1426}%
\special{pa 4892 1452}%
\special{pa 4910 1480}%
\special{pa 4930 1506}%
\special{pa 4950 1532}%
\special{pa 4972 1560}%
\special{pa 4994 1586}%
\special{pa 5018 1612}%
\special{pa 5042 1638}%
\special{pa 5068 1666}%
\special{pa 5094 1692}%
\special{pa 5120 1718}%
\special{pa 5146 1746}%
\special{pa 5172 1772}%
\special{pa 5198 1798}%
\special{pa 5200 1800}%
\special{sp}%
\special{pn 8}%
\special{sh 1.000}%
\special{ar 4000 600 50 50  0.0000000 6.2831853}%
\special{pn 8}%
\special{sh 1.000}%
\special{ar 5200 1800 50 50  0.0000000 6.2831853}%
\special{pn 8}%
\special{sh 0}%
\special{ar 5200 600 50 50  0.0000000 6.2831853}%
\special{pn 8}%
\special{sh 0}%
\special{ar 4800 1200 50 50  0.0000000 6.2831853}%
\special{pn 8}%
\special{sh 0}%
\special{ar 4400 1200 50 50  0.0000000 6.2831853}%
\special{pn 8}%
\special{sh 0}%
\special{ar 4000 1800 50 50  0.0000000 6.2831853}%
\put(25.0000,-6.5000){\makebox(0,0)[lb]{$V_n^{(-1)}$}}%
\put(25.0000,-18.5000){\makebox(0,0)[lb]{$V_n^{(1)}$}}%
\put(25.0000,-12.5000){\makebox(0,0)[lb]{$V_n^{(\alpha)}$}}%
\put(53.0000,-6.5000){\makebox(0,0)[lb]{$\cycmodp n1$}}%
\put(53.0000,-18.5000){\makebox(0,0)[lb]{$\cycmodp n{-1}$}}%
\put(49.0000,-12.5000){\makebox(0,0)[lb]{$\cycmodp n\alpha$}}%
\special{pn 8}%
\special{sh 1.000}%
\special{ar 3200 2200 50 50  0.0000000 6.2831853}%
\special{pn 8}%
\special{sh 0}%
\special{ar 3200 2400 50 50  0.0000000 6.2831853}%
\put(33.5000,-22.5000){\makebox(0,0)[lb]{irreducible}}%
\put(33.5000,-24.5000){\makebox(0,0)[lb]{reducible}}%
\put(20.0000,-4.0000){\makebox(0,0)[lb]{\fbox{Classical Case}}}%
\put(41.0000,-4.0000){\makebox(0,0)[lb]{\fbox{Quantum Case}}}%
\put(39.0000,-6.7000){\makebox(0,0)[rb]{$\cycmod n{-1}$}}%
\put(39.0000,-18.7000){\makebox(0,0)[rb]{$\cycmod n1$}}%
\put(43.0000,-12.7000){\makebox(0,0)[rb]{$\cycmod n\alpha$}}%
\special{pn 8}%
\special{pa 4450 1200}%
\special{pa 4750 1200}%
\special{dt 0.045}%
\put(46.0000,-11.3000){\makebox(0,0){$\cong$}}%
\put(18.0000,-6.5000){\makebox(0,0)[rb]{(skew-symmetric tensor)}}%
\put(18.0000,-18.5000){\makebox(0,0)[rb]{(symmetric tensor)}}%
\special{pn 8}%
\special{pa 4000 600}%
\special{pa 5200 1800}%
\special{da 0.070}%
\special{pn 8}%
\special{pa 2950 1200}%
\special{pa 3750 1200}%
\special{fp}%
\special{sh 1}%
\special{pa 3750 1200}%
\special{pa 3684 1180}%
\special{pa 3698 1200}%
\special{pa 3684 1220}%
\special{pa 3750 1200}%
\special{fp}%
\put(28.5000,-11.5000){\makebox(0,0)[lb]{{\footnotesize `quantum covering'}}}%
\put(18.0000,-12.5000){\makebox(0,0)[rb]{(generic $\alpha$)}}%
\end{picture}%
\medskip
\end{flushleft}
By looking at this `quantum covering structure' of interpolating property
as well as examining examples
(Examples 
\ref{ex:description_for_n=2_case}, 
\ref{ex:description_for_n=3_case} 
and
\ref{ex:per:description_for_n=2_case}, 
\ref{ex:per:description_for_n=3_case} 
respectively),
we infer an existence of `reciprocity' between cyclic modules
generated by $\Det$ and $\Per$.
Based on this observation,
we establish another conjecture (Conjecture B)
which describes a `reciprocity' between the multiplicities
of the irreducible summands of the cyclic modules
$\cycmod n\alpha$ and $\cycmodp n\alpha$.
In the very final position of the paper,
we introduce the partition functions
as respective generating functions
\begin{align*}
\Dtheta(t,\alpha)&\deq\sum_{\lambda:\text{partitions}}
\frac{\mult\lambda\alpha}{f^\lambda}t^{\abs\lambda}
=\sum_{n=0}^\infty\sum_{\lambda\vdash n}
\frac{\mult\lambda\alpha}{f^\lambda}t^n,\\
\Ptheta(t,\alpha)&\deq\sum_{\lambda:\text{partitions}}
\frac{\multp\lambda\alpha}{f^\lambda}t^{\abs\lambda}
=\sum_{n=0}^\infty\sum_{\lambda\vdash n}
\frac{\multp\lambda\alpha}{f^\lambda}t^n
\end{align*}
of the multiplicities
and restate certain weaker version
of two conjectures above in terms of the partition functions
(Here $\multp\lambda\alpha$ denotes the multiplicity of $\Umod n{\lambda}$
in $\cycmodp n\alpha$).
If $\alpha$ is not singular,
then $\Dtheta(t,\alpha)$ is identical
with the generating function $\prod_{n=1}^\infty(1-t^n)^{-1}$
of the number of partitions.

Since one can obtain neither the zeros of $q$-content discriminants
nor the discriminants themselves explicitly,
in order to understand a deeper structure of
the $q$-content discriminants as polynomials,
it might be an inevitable task to characterize
the $q$-content transition matrices in a suitable way,
e.g. like in a $R$-matrix formalism.
Although the situation is quite mysterious, in Proposition \ref{lem:symmetry},
we describe a symmetric structure of some matrix (see Lemma \ref{lem:reduction_1})
closely related to the $q$-content transition matrices.

\subsection*{Conventions}

Throughout the paper,
$\C$ is the complex number field and $\Z$ is the ring of rational integers.
The symbol $q$ denotes a nonzero complex number
and we always fix a branch for the square root $q^{\frac12}$.
We also assume that $q$ is not a root of unity
to assure the complete reducibility
of the finite dimensional representations of $\Uq$ and $\Hq$
as well as to utilize the highest weight theory of $\Uq$.
(We will further suppose that $q$ is `generic' to make the discussion simple.
See the end of Section \ref{QAD}.)

We denote by $\sym n$ the symmetric group of degree $n$.
The simple transposition $(k,k+1)$ is denoted by $s_k$.
For $w\in\sym n$, $\inv w$ is the inversion number of $w$
and $\cyc w$ is the number of cycles in $w$.
Notice that $\cyc\cdot$ is a class function on $\sym n$,
but $\inv\cdot$ is not.

For a partition (or a Young diagram)
$\lambda=(\lambda_1,\dots,\lambda_k)$ of $n$,
$\sym\lambda\deq\sym{\lambda_1}\times\dotsb\times\sym{\lambda_k}$
is the Young subgroup of $\sym n$.
We also denote by $\STab(\lambda)$ the set of
all standard tableaux with shape $\lambda$,
and put $f^\lambda=\card{\STab(\lambda)}$.

\section{Preliminaries on representations of quantum groups}
\label{Preliminaries}

We briefly recall the basic notion of
the quantum enveloping algebra $\Uq$,
the quantum matrix algebra $\Aq$,
and the Iwahori-Hecke algebra $\Hq$
to fix the conventions.
The conventions on quantum groups are the same
as in Noumi-Yamada-Mimachi \cite{NYM1993} and Jimbo \cite{J1986},
but the conventions for the Iwahori-Hecke algebra are
slightly different from those in Gyoja \cite{G1986}
because of the compatibility with the conventions on quantum algebras.

\subsection{Quantum enveloping algebra}

By definition, $\lattice$ is a $\Z$-module
$\lattice\deq\Z\e_1+\dots+\Z\e_n$
generated by the symbols $\e_1,\dots,\e_n$.
We fix a bilinear form $\blf{\cdot}{\cdot}$ on $\lattice$
defined by $\blf{\e_i}{\e_j}=\delta_{ij}$.
Each element in the lattice $\lattice$ is called an \emph{integral weight}.
The quantum enveloping algebra $\Uq$ is a $\C$-associative algebra
generated by the symbols $e_i,\, f_i$ ($1\le i\le n-1$)
and $q^\lambda$ ($\lambda\in\frac12\lattice$)
satisfying certain fundamental relations
(see \cite{NYM1993}).
The algebra $\Uq$ has a Hopf algebra structure
with the coproducts
\begin{align*}
&\cop(q^\lambda)=q^\lambda\otimes q^\lambda,\\
&\cop(e_k)=e_k\otimes q^{-(\e_k-\e_{k+1})/2}+q^{(\e_k-\e_{k+1})/2}\otimes e_k,\\
&\cop(f_k)=f_k\otimes q^{-(\e_k-\e_{k+1})/2}+q^{(\e_k-\e_{k+1})/2}\otimes f_k.
\end{align*}

The \emph{vector representation} $\rhoCn$ of $\Uq$ on
$\C^n$ is defined by
\begin{equation*}
\begin{split}
\rhoCn(q^\lambda)\cdot \ve_{j}=q^{\blf{\lambda}{\e_j}}\ve_{j},\qquad
\rhoCn(e_k)\cdot \ve_{j}=\delta_{j,k+1}\ve_{k},\qquad
\rhoCn(f_k)\cdot \ve_{j}=\delta_{jk}\ve_{k+1},
\end{split}
\end{equation*}
where $\{\ve_j\}_{1\le j\le n}$ is the standard basis of $\C^n$.
By the coproduct on $\Uq$,
the tensor product representation $\rhoCnn$
of $\Uq$ on $(\C^n)^{\otimes n}$
is given by
\begin{equation}\label{TensorAction}
\begin{split}
\rhoCnn(q^\lambda)\cdot \ve_{j_1}\otimes\dots\otimes \ve_{j_n}
&= q^{\blf{\lambda}{\e_{j_1}+\dots+\e_{j_n}}}\ve_{j_1}\otimes\dots\otimes \ve_{j_n},\\
\rhoCnn(e_k)\cdot \ve_{j_1}\otimes\dotsb\otimes \ve_{j_n}
&= \sum_{l=1}^n \delta_{j_l,k+1}q_k^l(j_1,\dots,j_n)
\ve_{j_1}\otimes\dotsb \ve_{j_{l-1}}\otimes \ve_{k}\otimes \ve_{j_{l+1}}
\otimes\dotsb\otimes \ve_{j_n},\\
\rhoCnn(f_k)\cdot \ve_{j_1}\otimes\dots\otimes \ve_{j_n}
&= \sum_{l=1}^n \delta_{j_l,k}q_k^l(j_1,\dots,j_n)
\ve_{j_1}\otimes\dotsb \ve_{j_{l-1}}\otimes \ve_{k+1}\otimes \ve_{j_{l+1}}
\otimes\dotsb\otimes \ve_{j_n},
\end{split}
\end{equation}
where we put $\q kl{j_1,\dots,j_n}\deq
q^{\blf{(\e_k-\e_{k+1})/2}{\e_{j_1}
+\dots+\e_{j_{l-1}}-\e_{j_{l+1}}-\dots-\e_{j_n}}}$  for simplicity.

Each finite dimensional irreducible $\Uq$-module
is a highest weight module,
and it is parametrized by a \emph{dominant} integral weight,
that is, an integral weight
$\lambda_1\e_1+\dotsb+\lambda_n\e_n$
with the property $\lambda_1\ge\dots\lambda_n\ge0$.
We identify a dominant integral weight $\lambda$
with a partition (or a Young diagram)
$(\lambda_1,\dots,\lambda_n)$.
We often denote by the same symbol $\lambda$
to indicate both the weight $\lambda_1\e_1+\dotsb+\lambda_n\e_n$
and the partition $(\lambda_1,\dots,\lambda_n)$.
The highest weight $\Uq$-module
corresponding to $\lambda$ is denoted by $\Umod n{\lambda}$.

\subsection{Quantum matrix algebra}

The quantum matrix algebra $\Aq$ is a $\C$-associative algebra
generated by $n^2$ letters $x_{ij}$ ($1\le i,j\le n$)
obeying the following fundamental relations
\begin{equation*}
\begin{split}
&x_{ik}x_{jk}=q x_{jk}x_{ik}, \quad
x_{ki}x_{kj}=q x_{kj}x_{ki} \qquad (i<j),\\
&x_{il}x_{jk}=x_{jk}x_{il},\quad
x_{ik}x_{jl}-x_{jl}x_{ik}=(q-q^{-1})x_{il}x_{jk} \qquad (i<j, k<l).
\end{split}
\end{equation*}
The algebra $\Aq$ becomes a bialgebra
having the coproduct
\begin{equation*}
\cop(x_{ij})=\sum_{k=1}^n x_{ik}\otimes x_{kj}.
\end{equation*}
The algebra $\Aq$ becomes a left $\Uq$-module by
\begin{equation*}
\begin{split}
\rho(q^\lambda)\cdot x_{ij}=q^{\blf{\lambda}{\e_j}}x_{ij},\qquad
\rho(e_k)\cdot x_{ij}=\delta_{j,k+1}x_{ik},\qquad
\rho(f_k)\cdot x_{ij}=\delta_{jk}x_{i,k+1}.
\end{split}
\end{equation*}
Using the coproduct of $\Uq$, (via the tensor product
representation) we have
\begin{equation}\label{eq:Uq-action_on_monomials}
\begin{split}
&\rho(q^\lambda)\cdot x_{i_1j_1}\dots x_{i_nj_n}
=q^{\blf{\lambda}{\e_{j_1}+\dots+\e_{j_n}}}x_{i_1j_1}\dots x_{i_nj_n},\\
&\rho(e_k)\cdot x_{i_1j_1}\dotsb x_{i_lj_l} \dotsb x_{i_nj_n}
=\sum_{l=1}^n \delta_{j_l,k+1}\cdot
\q kl{j_1,\dots,j_n}\cdot
x_{i_1j_1}\dotsb
{x_{i_lk}}
\dotsb x_{i_nj_n},\\
&\rho(f_k)\cdot x_{i_1j_1}\dotsb x_{i_lj_l} \dotsb x_{i_nj_n}
=\sum_{l=1}^n \delta_{j_l,k}\cdot
\q kl{j_1,\dots,j_n}\cdot
x_{i_1j_1}\dotsb
{x_{i_l,k+1}}
\dotsb x_{i_nj_n}.
\end{split}
\end{equation}
Notice that the $\Uq$-submodule $\bigoplus_{j=1}^n\C\cdot x_{ij}$
($i=1,\dots,n$) is equivalent to the vector representation $\C^n$.

\subsection{Iwahori-Hecke algebra and Schur-Weyl type duality}\label{IHA}

The Iwahori-Hecke algebra $\Hq$ is an associative $\C$-algebra
generated by the symbols $h_i$ ($1\le i\le n-1$)
with the fundamental relations
\begin{equation*}
\begin{split}
&h_ih_{i+1}h_i=h_{i+1}h_ih_{i+1}\quad(1\le i\le n-2),\\
&h_ih_j=h_jh_i\quad(\abs{i-j}\ge2),\\
&(h_i+q)(h_i-q^{-1})=0 \quad(1\le i\le n-1).
\end{split}
\end{equation*}
If $w=s_{i_1}\dotsb s_{i_l}$ is the shortest expression of $w\in\sym n$,
then we put $h_w=h_{i_1}\dotsb h_{i_l}$.
The elements $h_w$ for $w\in\sym n$ form a basis of $\Hq$
as a vector space.
The algebra $\Hq$ acts on $(\C^n)^{\otimes n}$ by
\begin{equation}\label{eq:action_of_Uq_on_Cnn}
\ve_{j_1}\otimes\dots\otimes \ve_{j_n}\cdot\pi(h_k)=
\begin{cases}
\ve_{j_1}\otimes\dots\otimes \ve_{j_{k+1}}\otimes \ve_{j_k}
\otimes\dots\otimes \ve_{j_n} & j_k<j_{k+1},\\
q^{-1}\ve_{j_1}\otimes\dots\otimes \ve_{j_n} & j_k=j_{k+1},\\
\ve_{j_1}\otimes\dots\otimes \ve_{j_{k+1}}\otimes \ve_{j_k}\otimes\dots\otimes \ve_{j_n}
-(q-q^{-1})\ve_{j_1}\otimes\dots\otimes \ve_{j_n} & j_k>j_{k+1}.
\end{cases}
\end{equation}
The subalgebra $\pi(\Hq)$ is the commutant of
$\rhoCnn(\Uq)$ in $\End((\C^n)^{\otimes n})$,
and vice versa (see, e.g. \cite{J1986}).
We have consequently the decomposition
\begin{equation*}
(\C^n)^{\otimes n}
\cong\bigoplus_{\lambda\vdash n}(\Umod n\lambda)^{\oplus f^\lambda}
\end{equation*}
as a $\Uq$-module.
This fact is referred to as \emph{Schur-Weyl duality}.

For a given Young diagram $\lambda\vdash n$,
we define
\begin{equation*}
e_+=e_+(\lambda)\deq\sum_{w\in W_+(\lambda)}q^{-\inv w}h_w,\qquad
e_-=e_-(\lambda)\deq\sum_{w\in W_-(\lambda)}(-q)^{\inv w}h_w,
\end{equation*}
where $W_{\pm}(\lambda)$ are certain subgroups of $\sym n$
(see \cite{G1986} for definition).
These satisfy the equations
\begin{equation*}
e_+^2=\kakko[3]{\sum_{w\in W_+(\lambda)}q^{-2\inv w}}e_+,\qquad
e_-^2=\kakko[3]{\sum_{w\in W_-(\lambda)}(-q)^{2\inv w}}e_-.
\end{equation*}
Using $e_\pm$,
we can define the $q$-Young symmetrizer $\yngsym(T)\in\Hq$
for each $T\in\STab(\lambda)$.
We refer to Gyoja \cite{G1986} for
precise and detailed information on
$q$-Young symmetrizers,
and we only give several examples here.

\begin{ex}
The $q$-Young symmetrizers for $3$-box standard tableaux
are given by
\begin{align*}
\yngsym(\minitab(123))&=e_+(\miniyng(3))
=1+q^{-1}h_1+q^{-1}h_2+q^{-3}h_1h_2h_1+q^{-2}h_1h_2+q^{-2}h_2h_1,\\
\yngsym(\minitab(12,3))&=h_2e_-(\miniyng(2,1))h_2^{-1}e_+(\miniyng(2,1))
=1+q^{-1}h_1-(q-q^{-1})h_2-q^{-1}h_1h_2h_1-h_1h_2-(1-q^{-2})h_2h_1,\\
\yngsym(\minitab(13,2))&=e_-(\miniyng(2,1))h_2e_+(\miniyng(2,1))h_2^{-1}
=1-qh_1-q^{-2}h_2h_1+q^{-1}h_1h_2h_1,\\
\yngsym(\minitab(1,2,3))&=e_-(\miniyng(1,1,1))
=1-qh_1-qh_2-q^{3}h_1h_2h_1+q^{2}h_1h_2+q^{2}h_2h_1.
\end{align*}
\eoe\end{ex}

\begin{ex}
For spaces of $q$-symmetric tensors and $q$-skew-symmetric tensors
representations, we have
\begin{align*}
\yngsym(\minitab(12{\cdot}n))=e_+((n))=\sum_{w\in\sym n}q^{-\inv w}h_w,\qquad
\yngsym(\minitab(1,2,\cdot,n))=e_-((1,\dots,1))=\sum_{w\in\sym n}(-q)^{\inv w}h_w
\end{align*}
in general.
\eoe\end{ex}

\section{Quantum $\alpha$-determinant}\label{QAD}

We now introduce a
\emph{quantum \pp{column} $\alpha$-determinant} $\Det$ by
\begin{equation*}
\Det\deq \sum_{w\in\sym{n}}
\alpha^{n-\cyc w}q^{\inv w}x_{w(1)1}\dotsb x_{w(n)n}
\in\Aq.
\end{equation*}
Since $(-1)^{n-\cyc w}=(-1)^{\inv w}$,
quantum $(-1)$-determinant
$\Det[-1]$ agrees with
the quantum determinant $\det_q$.
We also remark that $\Det(X)=\Det(\tX)$,
which follows from the fact that
$\cyc w=\cyc{w^{-1}}$,
$\inv w=\inv{w^{-1}}$ and
$x_{w(1)1}\dotsb x_{w(n)n}=x_{1w^{-1}(1)}\dotsb x_{nw^{-1}(n)}$
for any $w\in\sym n$.
For the sake of convenience, we write
\begin{equation*}
\qad{j_1,\dots,j_n}\deq
\sum_{w\in\sym{n}}\alpha^{n-\cyc w}q^{\inv w}
x_{w(1),j_1}\dots x_{w(n),j_n},
\end{equation*}
for $1\leq j_1,\ldots,j_n\leq n$.
We notice that $\Det=\qad{1,2,\dots,n}$.

\subsection{Quantum $\alpha$-determinant cyclic
modules $\cycmod n\alpha$}

We are interested in the $\Uq$-module
$\cycmod n\alpha\deq\rho(\Uq)\cdot\Det$.
The basic fact is that
every quantum $\alpha$-determinant
$\qad{j_1,\dots,j_n}$ is contained in $\cycmod n\alpha$
(Proposition \ref{lem:q-analogue_of_Lemma_2.2}).
To prove this,
we first notice the following
\begin{lem}\label{lem:action_on_dets}
The equalities
\begin{equation*}
\begin{split}
\rho(q^\lambda)\cdot\qad{j_1,\dots,j_n}
&= q^{\blf{\lambda}{\e_{j_1}+\dots+\e_{j_n}}}\qad{j_1,\dots,j_n},\\
\rho(e_k)\cdot\qad{j_1,\dots,j_n}
&= \sum_{l=1}^n \delta_{j_l,k+1}q_k^l(j_1,\dots,j_n)
\qad{j_1,\dots,j_{l-1},k,j_{l+1},\dots,j_n},\\
\rho(f_k)\cdot\qad{j_1,\dots,j_n}
&= \sum_{l=1}^n \delta_{j_l,k}q_k^l(j_1,\dots,j_n)
\qad{j_1,\dots,j_{l-1},k+1,j_{l+1},\dots,j_n}
\end{split}
\end{equation*}
hold.
\end{lem}

\begin{proof}
By \eqref{eq:Uq-action_on_monomials},
the assertion is verified by a straightforward calculation.
\end{proof}

We give some instructive example which may indicate
a highest weight vector of a irreducible representation (see Section \ref{QAD}).
\begin{ex}\label{ex:hwv_for_(2,1)}
We have
\begin{equation*}
\begin{split}
\rho(e_1)\cdot\qad{1,1,2}
&=q_1^3(1,1,2)\qad{1,1,1}
=q^{\blf{(\e_1-\e_2)/2}{\e_1+\e_1}}\qad{1,1,1}=q\qad{1,1,1},\\
\rho(e_1)\cdot\qad{1,2,1}
&=q_1^2(1,2,1)\qad{1,1,1}
=q^{\blf{(\e_1-\e_2)/2}{\e_1-\e_1}}\qad{1,1,1}=\qad{1,1,1},\\
\rho(e_1)\cdot\qad{2,1,1}
&=q_1^1(2,1,1)\qad{1,1,1}
=q^{\blf{(\e_1-\e_2)/2}{-\e_1-\e_1}}\qad{1,1,1}=q^{-1}\qad{1,1,1},
\end{split}
\end{equation*}
and hence, we conclude that
\begin{equation*}
\rho(e_1)\cdot(\qad{1,1,2}-q\qad{1,2,1})
=\rho(e_1)\cdot(\qad{1,2,1}-q\qad{2,1,1})
=0.
\end{equation*}
These two vectors are also killed by $\rho(e_2)$ trivially.
\eoe\end{ex}

\begin{prop}[Quantum analogue of {\cite[Lemma 2.2]{MW2005}}]%
\label{lem:q-analogue_of_Lemma_2.2}
The equality
\begin{equation}\label{eq:V_as_linear_span}
\cycmod n\alpha=\sum_{1\le j_1,\dots,j_n \le n}\C\cdot\qad{j_1,\dots,j_n}
\end{equation}
holds.
\end{prop}

\begin{proof}
Let $\linspan n\alpha$ be the
right-hand side of \eqref{eq:V_as_linear_span}.
By Lemma \ref{lem:action_on_dets},
$\linspan n\alpha$ is $\rho(\Uq)$-invariant and
$\cycmod n\alpha\subset \linspan n\alpha$.
To prove the opposite inclusion
$\cycmod n\alpha\supset \linspan n\alpha$,
we introduce the linear map
\begin{equation}\label{eq:def_of_intertwiner}
\Phi_{n,q}^{(\alpha)}:
(\C^n)^{\otimes n}\ni
\ve_{j_1}\otimes\dotsb\otimes \ve_{j_n}
\longmapsto
\qad{j_1,\dots,j_n}
\in\linspan n\alpha
\qquad(1\le j_1,\dots,j_n\le n).
\end{equation}
By Lemma \ref{lem:action_on_dets} again and
the formula \eqref{TensorAction}, 
we find that $\Phi_{n,q}^{(\alpha)}$ defines a surjective $\Uq$-intertwiner
such that
$\Phi_{n,q}^{(\alpha)}(\ve_1\otimes\dots\otimes\ve_n)=\Det$.
Using the elementary fact that
\begin{equation*}
(\C^n)^{\otimes n}=\rhoCnn(\Uq)\cdot\ve_1\otimes\dots\otimes\ve_n,
\end{equation*}
we have
\begin{equation*}
\cycmod n\alpha=\rho(\Uq)\cdot\Det
=\Phi_{n,q}^{(\alpha)}(\rhoCnn(\Uq)\cdot\ve_1\otimes\dots\otimes\ve_n)
=\Phi_{n,q}^{(\alpha)}((\C^n)^{\otimes n})
\supset\linspan n\alpha.
\end{equation*}
This completes the proof.
\end{proof}

By Proposition \ref{lem:q-analogue_of_Lemma_2.2}
and the surjectivity of the intertwiner $\Phi_{n,q}^{(\alpha)}$
defined in \eqref{eq:def_of_intertwiner},
the cyclic $\Uq$-module $\cycmod n\alpha$
is isomorphic to the tensor product module $(\C^n)^{\otimes n}$
if and only if the intertwiner $\Phi_{n,q}^{(\alpha)}$ is bijective,
that is, the $\alpha$-determinants $\qad{j_1,\dots,j_n}$ are
linearly independent.
Namely, we have the following basic result.
\begin{prop}\label{prop:ird_generic}
Put
\begin{equation*}
\singular n\deq
\set[1]{\alpha\in\C}
{\text{\upshape$\qad{j_1,\dots,j_n}$ are linearly dependent}}.
\end{equation*}
If $\alpha\in\C\setminus\singular n$,
then the irreducible decomposition of
$\cycmod n\alpha$ is given as
\begin{equation*}
\cycmod n\alpha
\cong\bigoplus_{\lambda\vdash n}(\Umod n\lambda)^{\oplus f^\lambda}.
\end{equation*}
In other words, the \emph{multiplicity} $\mult\lambda\alpha$
of the irreducible representation
$\Umod n\lambda$ in $\cycmod n\alpha$
is equal to $f^\lambda$.
\qed
\end{prop}

Let us look at the weight space decomposition of $\cycmod n\alpha$.
For each $\vj=(j_1,\dots,j_n)\in\{1,\dots,n\}^n$,
we associate an integral weight
\begin{equation*}
\wt(\vj)=\sum_{i=1}^n \e_{j_i}=
\sum_{i=1}^n \nu_i\e_i
\in\lattice
\qquad\bigl(\nu_k=\card{\set{l}{j_l=k}}\bigr).
\end{equation*}
Lemma \ref{lem:action_on_dets} says that
$\qad{\vj}=\qad{j_1,\dots,j_n}$ is a weight vector of weight $\wt(\vj)$.
For convenience, we put
\begin{align*}
\ulattice&\deq
\set[1]{\nu=\nu_1\e_1+\dotsb+\nu_n\e_n\in\lattice}
{\nu_j\ge0,\,\nu_1+\dotsb+\nu_n=n},\\
I_n(\lambda)&\deq
\set[1]{\vi\in\{1,\dots,n\}^n}{\wt(\vi)=\lambda}
\qquad(\lambda\in\ulattice).
\end{align*}
By Lemma \ref{lem:action_on_dets}, we have the
\begin{lem}\label{lem:weight_decomp}
For an integral weight $\nu\in\ulattice$,
define the subspace
\begin{align*}
\cycmod n\alpha(\nu)&\deq
\sum_{\vi\in I_n(\nu)}\C\cdot\qad{\vi}
\end{align*}
consisting of all weight vectors of weight $\nu$.
Then the following decomposition
\begin{equation*}
\cycmod n\alpha=
\bigoplus_{\nu\in\ulattice}\cycmod n\alpha(\nu)
\end{equation*}
holds. \qed
\end{lem}

\subsection{Singular points for the decomposition}

By Proposition \ref{prop:ird_generic},
what remains important is
to study the set $\singular n$
and determine the irreducible decomposition of $\cycmod n\alpha$
for $\alpha\in\singular n$.
We first give another description of $\singular n$
as a set of zeros of a certain polynomial defined below.

We notice that
$\qad[0]{j_1,\dots,j_n}$
is nothing but the monomial
$x_{1j_1}\dotsb x_{nj_n}$.
It hence follows that
the vectors $\qad[0]{j_1,\dots,j_n}$ are
linearly independent
(i.e. $0\notin\singular n$)
and each $\alpha$-determinant $\qad{j_1,\dots,j_n}$ is
a linear combination of the monomials
$\qad[0]{j_1,\dots,j_n}$, say
\begin{equation*}
\qad{j_1,\dots,j_n}=\sum_{1\le i_1,\dots,i_n \le n}
\uF_{n,q}(\alpha;{i_1,\dots,i_n;j_1,\dots,j_n})\qad[0]{i_1,\dots,i_n}
\end{equation*}
for some $\uF_{n,q}(\alpha;{i_1,\dots,i_n;j_1,\dots,j_n})\in\C$.
It is immediate to see that each
$\uF_{n,q}(\alpha;{\vi;\vj})$ ($\vi,\vj\in\{1,2,\dots,n\}^n$)
is a polynomial in $\alpha$ and $q$
with integral coefficients.

Consider the $n^n\times n^n$ matrix
$\uF_{n,q}(\alpha)
\deq(\uF_{n,q}(\alpha;{\vi;\vj}))_{\vi,\vj\in\{1,\dots,n\}^n}$.
The determinant $\uC_{n,q}(\alpha)\deq\det\uF_{n,q}(\alpha)$
is a polynomial in $\alpha$ and $q$
with integral coefficients,
and it is not identically zero because $\uF_{n,q}(0)$ is the identity matrix.
Thus we have the
\begin{lem}
The set $\singular n$ is given by
\begin{equation*}
\singular n=\set[1]{\alpha\in\C}{\uC_{n,q}(\alpha)=0}.
\end{equation*}
In particular,
$\singular n$ is a finite set.
\qed
\end{lem}

The cardinality $\card{\singular n}$ does depend on the parameter $q$.
In what follows, for simplicity,
we further impose an assumption on the parameter $q$ that
the value $q$ maximize $\card{\singular n}$
as a function in $q$.
(It is sufficient to assume that $q$ is transcendental, for instance.)

Let us put
\begin{equation*}
\uF_{q}^\lambda(\alpha)
\deq\Bigl(\uF_{n,q}(\alpha;{\vi;\vj})\Bigr)_{\vi,\vj\in I_n(\lambda)}
\end{equation*}
for an integral weight $\lambda\in\ulattice$.
Then, by Lemma \ref{lem:weight_decomp},
the matrix $\uF_{n,q}(\alpha)$ is a direct sum
\begin{equation*}
\uF_{n,q}(\alpha)\sim\bigoplus_{\lambda\in\ulattice}\uF_{q}^\lambda(\alpha)
\end{equation*}
of the smaller matrices $\uF_{q}^\lambda(\alpha)$
because $\uF_{n,q}(\alpha;{\vi;\vj})=0$ if $\wt(\vi)\ne\wt(\vj)$.
Here, for given square matrices $A$ and $B$,
we write $A\sim B$ when $B=PAP^{-1}$ for some invertible matrix $P$.
If we put $\uC_{q}^\lambda(\alpha)=\det\uF_{q}^\lambda(\alpha)$,
it is clear that
\begin{equation}\label{eq:uC_decomp}
\uC_{n,q}(\alpha)=\prod_{\lambda\in\ulattice}\uC_{q}^\lambda(\alpha).
\end{equation}

The following lemma is immediately verified.
\begin{lem}\label{lem:reduced_to_dominant_weight}
If $\blf{\lambda}{\e_i}=\blf{\mu}{\e_{\sigma(i)}}$
\pp{$1\le i\le n$}
for some $\sigma\in\sym n$, then
$\uF_{q}^\lambda(\alpha)\sim\uF_{q}^\mu(\alpha)$
\pp{$\lambda,\,\mu\in\ulattice$}.
In particular,
for each $\mu\in\ulattice$,
there exists a unique dominant integral weight
$\lambda\in\ulattice$
such that $\uF_{q}^\lambda(\alpha)\sim\uF_{q}^\mu(\alpha)$.
\qed
\end{lem}
Consequently,
together with \eqref{eq:uC_decomp},
we have
\begin{equation*}
\singular n
=\bigcup_{\lambda\in\ulattice}
\set[1]{\alpha\in\C}{\uC_{q}^\lambda(\alpha)=0}
=\bigcup_{\lambda\in\domlattice}
\set[1]{\alpha\in\C}{\uC_{q}^\lambda(\alpha)=0}.
\end{equation*}
Here $\domlattice$ is the set of dominant weights in $\ulattice$.
We will regard a dominant integral weight
$\lambda=\lambda_1\e_1+\dots+\lambda_n\e_n\in\domlattice$
as a partition (or a Young diagram)
$(\lambda_1,\dots,\lambda_n)\vdash n$.
Also, we sometimes write $\lambda\vdash n$
to indicate $\lambda\in\domlattice$.

It seems quite difficult in general to
determine the polynomials $\uC_q^\lambda(\alpha)$
as well as the matrix $\uF_{q}^\lambda(\alpha)$
explicitly.
Therefore,
any characterization of the matrix $\uF_{q}^\lambda(\alpha)$,
for instance, either by difference equations or
in the framework of $R$-matrices would be interesting (if any).
From this point of view, the following property of
the matrices $\uF_{q}^\lambda(\alpha)$ is considerably remarkable.
(See Examples \ref{ex:(2,1)}, and \ref{ex:(3,1)} in \S4.6).
\begin{prop}\label{lem:symmetry}
The matrix $\uF_{q}^\lambda(\alpha)$ is symmetric.
\end{prop}
Before proceeding to the proof,
we prepare several convention.
We associate a sequence
\begin{equation*}
\vk(\lambda)=(k_1(\lambda),\dots,k_n(\lambda))\deq
(\overbrace{1,\dots,1}^{\lambda_1},
\overbrace{2,\dots,2}^{\lambda_2},\dots,
\overbrace{n,\dots,n}^{\lambda_n})\in I_n(\lambda)
\end{equation*}
to each integral weight
$\lambda=\lambda_1\e_1+\lambda_2\e_2+\dots+\lambda_n\e_n\in\ulattice$.
We define a right $\sym n$-action on $I_n(\lambda)$ by
$\vi^\sigma=(i_{\sigma(1)},\dots,i_{\sigma(n)})$
for $\sigma\in\sym n$ and $\vi=(i_1,\dots,i_n)\in I_n(\lambda)$.
Notice that the right action $I_n(\lambda)\tca\sym n$ is transitive and
the stabilizer of $\vk(\lambda)$
is the Young subgroup $\sym\lambda$.
Therefore we have
\begin{equation*}
I_n(\lambda)\cong \sym\lambda\backslash\sym n
\quad\text{and/or}\quad
I_n(\lambda)=\vk(\lambda)\cdot\sym n.
\end{equation*}
Thus we have another expression
\begin{equation*}
\uF_{q}^\lambda(\alpha)\sim
\Bigl(\uF_q^\lambda(\alpha;{\tau,\sigma})\Bigr)%
_{\tau,\sigma\in\sym\lambda\backslash\sym n},
\end{equation*}
where we put
$\uF_q^\lambda(\alpha;{\tau,\sigma})
\deq \uF_{n,q}(\alpha;{\vk(\lambda)^\tau;\vk(\lambda)^\sigma})$.

\begin{proof}[Proof of Proposition \ref{lem:symmetry}]
For convenience, we put
\begin{equation*}
X_{\vk}(g,\sigma)\deq x_{g(1)k_{\sigma(1)}}\dotsb x_{g(n)k_{\sigma(n)}}
\end{equation*}
for $g\in\sym n$ and $\sigma\in\sym\lambda\backslash\sym n$.
We also define $f^g_{\tau,\sigma}(\lambda)$ by
\begin{equation*}
q^{\inv g}X_{\vk}(g,\sigma)=
\sum_{\tau\in\sym\lambda\backslash\sym n}
f^g_{\tau,\sigma}(\lambda) X_{\vk}(1,\tau)
\end{equation*}
for $g\in\sym n$ and $\sigma,\tau\in\sym\lambda\backslash\sym n$.
It follows that
\begin{equation*}
\uF_q^\lambda(\alpha;{\tau,\sigma})
=\sum_{g\in\sym n}\alpha^{n-\cyc g}f^g_{\tau,\sigma}(\lambda).
\end{equation*}
Suppose that a permutation $g\in\sym n$ and a simple transposition $s_i$
satisfies the condition $\inv{gs_i}>\inv g$,
which is equivalent to the condition $g(i)<g(i+1)$.
It follows that
\begin{equation*}
\begin{split}
X_{\vk}(gs_i,\sigma)
&=\begin{cases}
q^{-1}X_{\vk}(g,\sigma s_i) & k_{\sigma(i)}=k_{\sigma(i+1)},\\
X_{\vk}(g,\sigma s_i) & k_{\sigma(i)}<k_{\sigma(i+1)},\\
X_{\vk}(g,\sigma s_i)-(q-q^{-1})X_{\vk}(g,\sigma) & k_{\sigma(i)}>k_{\sigma(i+1)}.
\end{cases}
\end{split}
\end{equation*}
This yields the relation
\begin{equation*}
\begin{split}
f^{gs_i}_{\tau,\sigma}(\lambda)
=\theta_0(\sigma,i)f^g_{\tau,\sigma}(\lambda)
+\theta_1(\sigma,i)f^g_{\tau,\sigma s_i}(\lambda)
\end{split}
\end{equation*}
where we put
\begin{equation*}
\theta_0(w,i)\deq\begin{cases}
1-q^2 & k_{w(i)}>k_{w(i+1)},\\
0 & \text{otherwise},
\end{cases}
\qquad
\theta_1(w,i)\deq\begin{cases}
1 & k_{w(i)}=k_{w(i+1)},\\
q & \text{otherwise}.
\end{cases}
\end{equation*}
Therefore,
for a given permutation $g=s_{i_l}\dotsb s_{i_1}$ of length $l$,
we have
\begin{equation*}
f^{s_{i_l}\dotsb s_{i_1}}_{\tau,\sigma}(\lambda)
=\sum_{(j_1,\dots,j_l)\in\{0,1\}^l}
\Theta_{\vk}\!\left(\begin{matrix}
{i_1},\dots,{i_l}\\
j_1,\dots,j_l
\end{matrix}
;\sigma\right)
\delta_{\tau,\sigma s_{i_1}^{j_1}\dotsb s_{i_l}^{j_l}}^{\lambda}
\end{equation*}
by induction.
Here we put
\begin{equation*}
\begin{split}
\Theta_{\vk}\!\left(\begin{matrix}
{i_1},\dots,{i_l}\\
j_1,\dots,j_l
\end{matrix}
;\sigma\right)&\deq
\theta_{j_1}(\sigma,i_1)\theta_{j_2}(\sigma s_{i_1}^{j_1},i_2)%
\dotsb%
\theta_{j_l}(\sigma s_{i_1}^{j_1}\dotsb s_{i_{l-1}}^{j_{l-1}},i_l),\\
\delta_{\tau,\sigma}^{\lambda}&\deq
\begin{cases}
1 & \tau\sigma^{-1}\in\sym\lambda,\\
0 & \text{otherwise}.
\end{cases}
\end{split}
\end{equation*}
We notice that
\begin{equation*}
\theta_{j_p}(\sigma s_{i_1}^{j_1}\dotsb s_{i_{p-1}}^{j_{p-1}},i_p)
=\theta_{j_p}(\sigma s_{i_1}^{j_1}\dotsb s_{i_{p}}^{j_{p}},i_p)
\end{equation*}
for each $p=1,2,\dots,l$.
It therefore follows that
if $\tau^{-1}\sigma s_{i_1}^{j_1}\dotsb s_{i_l}^{j_l}\in\sym\lambda$,
then
\begin{equation*}
\begin{split}
\Theta_{\vk}\!\left(\begin{matrix}
{i_1},\dots,{i_l}\\
j_1,\dots,j_l
\end{matrix}
;\sigma\right)
&=\theta_{j_1}(\sigma,i_1)\theta_{j_2}(\sigma s_{i_1}^{j_1},i_2)%
\dotsb%
\theta_{j_l}(\sigma s_{i_1}^{j_1}\dotsb s_{i_{l-1}}^{j_{l-1}},i_l)\\
&=\theta_{j_1}(\sigma s_1^{j_1},i_1)
\theta_{j_2}(\sigma s_{i_1}^{j_1}s_{i_2}^{j_2},i_2)%
\dotsb%
\theta_{j_l}(\sigma s_{i_1}^{j_1}\dotsb s_{i_{l}}^{j_{l}},i_l)\\
&=\theta_{j_l}(\tau,i_l)\theta_{j_{l-1}}(\tau s_{i_l}^{j_l},i_{l-1})%
\dotsb%
\theta_{j_1}(\tau s_{i_l}^{j_l}\dotsb s_{i_{2}}^{j_{2}},i_1)
=\Theta_{\vk}\!\left(\begin{matrix}
{i_l},\dots,{i_1}\\
j_l,\dots,j_1
\end{matrix}
;\tau\right).
\end{split}
\end{equation*}
This immediately implies that
$f^g_{\tau,\sigma}(\lambda)=f^{g^{-1}}_{\sigma,\tau}(\lambda)$,
and hence the symmetry $\uF^\lambda_{\tau,\sigma}=\uF^\lambda_{\sigma,\tau}$ follows
as we desired.
\end{proof}

\subsection{Highest weight vectors in $\cycmod n\alpha$}

The aim of the present subsection is to construct a set of vectors
$\set[1]{\hwv T}{T\in\STab(\lambda),\ \lambda\vdash n}$
in $\cycmod n\alpha$ satisfying the following conditions:
(a) If $T\in\STab(\lambda)$,
then $\hwv T\in\cycmod n\alpha(\lambda)$,
(b) each $\hwv T$ is killed by $\rho(e_k)$ ($1\le k<n$),
(c) $\cycmod n\alpha=\bigoplus_{\lambda\vdash n}
\left\{\sum_{T\in\STab(\lambda)}
\rho(\Uq)\cdot \hwv T\right\}$.
To achieve this, we first construct such vectors $\hwv T$
for $\alpha\in\C\setminus\singular n$,
and then extend the definition of them to any $\alpha\in\C$.
So, we suppose that $\alpha\in\C\setminus\singular n$ for a while.

For a quantum $\alpha$-determinant $\qad{j_1,\dots,j_n}$
and $h_k\in\Hq$, define
\begin{equation*}
\begin{split}
\qad{j_1,\dots,j_n}\cdot\pi^{(\alpha)}(h_k)&=
\begin{cases}
\qad{j_1,\dots,j_{k+1},j_k,\dots,j_n} & j_k<j_{k+1},\\
q^{-1}\qad{j_1,\dots,j_n} & j_k=j_{k+1},\\
\qad{j_1,\dots,j_{k+1},j_k,\dots,j_n}
-(q-q^{-1})\qad{j_1,\dots,j_n} & j_k>j_{k+1}.
\end{cases}
\end{split}
\end{equation*}
Each $\pi^{(\alpha)}(h_k)$ is extended
as a linear operator on $\cycmod n\alpha$
and defines a right $\Hq$-module structure on $\cycmod n\alpha$.

\begin{rem}
When $\alpha\in\singular n$,
we cannot extend $\pi^{(\alpha)}(h_k)$
to a linear operator on $\cycmod n\alpha$
as we see in the following example:
When $\alpha=\frac1{q^3+q^2-q}\in\singular 3$,
we have a nontrivial linear relation
\begin{equation*}
\qad{1,1,2}+(1-q)\qad{1,2,1}-q\qad{2,1,1}=0.
\end{equation*}
However, since
\begin{align*}
\qad{1,1,2}\cdot\pi^{(\alpha)}(h_1)&=q^{-1}\qad{1,1,2},\qquad
\qad{1,2,1}\cdot\pi^{(\alpha)}(h_1)=\qad{2,1,1},\\
\qad{2,1,1}\cdot\pi^{(\alpha)}(h_1)&=\qad{1,2,1}-(q-q^{-1})\qad{2,1,1},
\end{align*}
it follows that
\begin{equation*}
\begin{split}
&\qad{1,1,2}\cdot\pi^{(\alpha)}(h_1)+(1-q)\qad{1,2,1}\cdot\pi^{(\alpha)}(h_1)
-q\qad{2,1,1}\cdot\pi^{(\alpha)}(h_1)\\
=&\frac{(1-q^2)(1-q+q^2)}{q(1-q-q^2)}
\Bigl(\qad[0]{1,1,2}+(1-q)\qad[0]{1,2,1}-q\qad[0]{2,1,1}\Bigr)\ne0,
\end{split}
\end{equation*}
which means that $\pi^{(\alpha)}(h_1)$
cannot be extended to a linear operator on $\cycmod3\alpha$
when $\alpha=\frac1{q^3+q^2-q}\in\singular 3$.
\end{rem}

It is directly checked that
\begin{equation*}
\Phi_{n,q}^{(\alpha)}(\ve_{j_1}\otimes\dotsb\otimes\ve_{j_n}\cdot\pi(h_k))
=\qad{j_1,\dots,j_n}\cdot\pi^{(\alpha)}(h_k).
\end{equation*}
Hence, each operator $\pi^{(\alpha)}(h_k)$ commutes with
the $\rho(\Uq)$-action.
In particular,
$\pi^{(\alpha)}(\Hq)$ is the commutant of $\rho(\Uq)$
in $\End \cycmod n\alpha$ and vice versa.

For a standard tableau $T\in\STab(\lambda)$ of size $n$ ($\lambda\vdash n$),
we define $\vj(T)=(j_1(T),\dots,j_n(T))\in I_n(\lambda)$ by
\begin{equation}\label{eq:JT}
j_p(T)=i \iff \text{the number written in the $(i,j)$-box in $T$ is $p$.}
\end{equation}
We set
\begin{equation*}
\hwv T\deq \qad{\vj(T)}\cdot\pi^{(\alpha)}(\yngsym(T))
\in\cycmod n\alpha(\lambda)
\qquad(T\in\STab(\lambda)),
\end{equation*}
where $\yngsym(T)$ is the $q$-Young symmetrizer for $T$
(see Section \ref{IHA}).
This is a highest weight vector of weight $\lambda$.
By definition, each vector $\hwv T$ has an expression
\begin{equation}\label{eq:vTa}
\hwv T=\sum_{\sigma\in\sym\lambda\backslash\sym n}
Q_T^\sigma(q)\qad{\vk(\lambda)^\sigma}
=\sum_{\tau\in\sym\lambda\backslash\sym n}
\left\{
\sum_{\sigma\in\sym\lambda\backslash\sym n}
Q_T^\sigma(q)\uF_q^\lambda(\alpha;{\tau,\sigma})
\right\}\qad[0]{\vk(\lambda)^\tau}
\end{equation}
for certain \emph{polynomials $Q_T^\sigma(q)$ in $q$}.
For later use, we define the
$f^\lambda\times\card{\sym\lambda\backslash\sym n}$ matrix $\uQ_n^\lambda(q)$ by
\begin{equation}\label{eq:matrix_Q}
\uQ_n^\lambda(q)
=\bigl(Q_T^\sigma(q)\bigr)_{T\in\STab(\lambda),\sigma\in\sym\lambda\backslash\sym n}.
\end{equation}

Similar to the classical case,
the vectors $\hwv T$ for $T\in\STab(\lambda)$ form a basis of the subspace
\begin{equation*}
W_{n,q}^{(\alpha)}(\lambda)\deq
\set[1]{v\in \cycmod n\alpha(\lambda)}
{\rho(e_k)\cdot v=0\ (1\le k<n)}
\end{equation*}
consisting of the highest weight vectors of highest weight $\lambda$.
Therefore,
the cyclic module $\rho(\Uq)\cdot \hwv T$ is
equivalent to $\Umod n\lambda$ for each $T\in\STab(\lambda)$
and we have
\begin{equation*}
\cycmod n\alpha=\bigoplus_{\lambda\vdash n}\bigoplus_{T\in\STab(\lambda)}
\rho(\Uq)\cdot \hwv T.
\end{equation*}
In particular, every quantum $\alpha$-determinant $\qad{i_1,\dots,i_n}$
is written in the form
\begin{equation}\label{eq:generate_the_module}
\qad{i_1,\dots,i_n}=\sum_{\lambda\vdash n}\sum_{T\in\STab(\lambda)}
\rho(a^{(\alpha)}(T))\cdot \hwv T
\qquad(\exists a^{(\alpha)}(T)\in\Uq[\alpha]).
\end{equation}
Here we notice that the right-hand side of \eqref{eq:vTa} makes sense
even if $\alpha\in\singular n$,
though the vector $\hwv T$ is defined only
for $\alpha\in\C\setminus\singular n$.
Actually, it is a linear combination of monomials $x_{1i_1}\dotsb x_{ni_n}$
whose coefficient is a polynomial in $\alpha$.
So we \emph{extend} the definition of $\hwv T$ for any $\alpha\in\C$
by the expression \eqref{eq:vTa}.

\begin{lem}\label{lem:highest weight}
The vector $\hwv T$ is a highest weight vector
in $\cycmod n\alpha(\lambda)$
whenever $\hwv T\ne0$.
\end{lem}
\begin{proof}
By definition, each vector $\rho(e_k)\cdot \hwv T$
is a polynomial function in $\alpha$.
Therefore, the property $\rho(e_k)\cdot \hwv T=0$ is equivalent
to some algebraic equation on $\alpha$.
Since any complex number in $\C\setminus\singular n$
is a root of the equation,
we have $\rho(e_k)\cdot \hwv T=0$ for any $\alpha\in\C$.
This completes the proof.
\end{proof}

The formula \eqref{eq:generate_the_module} is valid for all $\alpha\in\C$
by a similar `polynomial' discussion.
It hence follows that
\begin{equation*}
\cycmod n\alpha=\bigoplus_{\lambda\vdash n}
\ckakko[4]{\sum_{T\in\STab(\lambda)}
\rho(\Uq)\cdot \hwv T}
\end{equation*}
for \emph{any} $\alpha\in\C$.
We notice that
$\{\hwv T\}_{T\in\STab(\lambda)}$
\emph{generates} $W_{n,q}^{(\alpha)}(\lambda)$
and $\mult\lambda\alpha=\dim W_{n,q}^{(\alpha)}(\lambda)$.

\subsection{Explicit decomposition of
$\cycmod n\alpha$ --- examples for $n=2,3$}

\begin{ex}\label{ex:description_for_n=2_case}
Let us see the simplest case,
the $\Uq[2]$-module $\cycmod2\alpha$.
Since
\begin{align*}
&\qad{1,1}=(1+\alpha)\qad[0]{1,1},&&
\qad{1,2}=\qad[0]{1,2}+\alpha q\qad[0]{2,1},\\
&\qad{2,1}=\alpha q\qad[0]{1,2}+(1+\alpha-\alpha q^2)\qad[0]{2,1},&&
\qad{2,2}=(1+\alpha)\qad[0]{2,2},
\end{align*}
we have
\begin{equation*}
\uF_{2,q}(\alpha)=\begin{pmatrix}
1+\alpha & 0 & 0 & 0 \\
0 & 1 & \alpha q & 0 \\
0 & \alpha q & 1+\alpha-\alpha q^2 & 0 \\
0 & 0 & 0 & 1+\alpha
\end{pmatrix}
=\uF_{q}^{(2,0)}(\alpha)\oplus\uF_{q}^{(1,1)}(\alpha)
\oplus\uF_{q}^{(0,2)}(\alpha),
\end{equation*}
$\uC_{2,q}(\alpha)=\det\uF_{2,q}(\alpha)=(1+\alpha)^3(1-\alpha q^2)$
and $\singular 2=\{-1,q^{-2}\}$.

If $\alpha\in\C\setminus\singular 2$,
then we have $\cycmod2\alpha\cong \Umod2{(2)}\oplus \Umod2{(1,1)}$
by Proposition \ref{prop:ird_generic}.
Each irreducible component is explicitly written as
\begin{equation*}
\begin{split}
&\Umod2{(2)}:
0 \Larrow{f_1} \C\qad{2,2} \LRarrow{e_1}{f_1} \C(q\qad{1,2}+\qad{2,1})
\LRarrow{e_1}{f_1} \C\qad{1,1} \Rarrow{e_1} 0,\\
&\Umod2{(1,1)}:
0 \Larrow{f_1} \C(\qad{1,2}-q\qad{2,1}) \Rarrow{e_1} 0.
\end{split}
\end{equation*}
The highest weight vectors of these modules are
\begin{equation*}
\qad{1,1}=(1+\alpha)x_{11}x_{21},\qquad
\qad{1,2}-q\qad{2,1}=(1-\alpha q^2)\det_q.
\end{equation*}
Hence the component $\Umod2{(2)}$ (resp. $\Umod2{(1,1)}$) disappears
if $\alpha=-1$ (resp. $\alpha=q^{-2}$).
We have thus
\begin{equation*}
\cycmod2\alpha \cong \begin{cases}
\Umod2{(2)} & \alpha=-1,\\
\Umod2{(1,1)} & \alpha=q^{-2},\\
\Umod2{(2)}\oplus \Umod2{(1,1)} & \alpha\neq-1,q^{-2}.
\end{cases}
\end{equation*}
\eoe\end{ex}

\begin{ex}\label{ex:description_for_n=3_case}
Look at the $\Uq[3]$-module $\cycmod3\alpha$.
If we put
\begin{align*}
v^{(3)}=&\qad{1,1,1}
=(1+\alpha)(1+2\alpha)x_{11}x_{21}x_{31},\\
v_1^{(2,1)}=&\qad{1,1,2}+(1-q)\qad{1,2,1}-q\qad{2,1,1}\\
=&(1+\alpha)(1+(q-q^2-q^3)\alpha)
(x_{11}x_{21}x_{32}+(1-q)x_{11}x_{22}x_{31}-qx_{12}x_{21}x_{31}),\\
v_2^{(2,1)}=&\qad{1,1,2}-(1+q)\qad{1,2,1}+q\qad{2,1,1}\\
=&(1+\alpha)(1+(-q-q^2+q^3)\alpha)
(x_{11}x_{21}x_{32}-(1+q)x_{11}x_{22}x_{31}+qx_{12}x_{21}x_{31}),\\
v^{(1,1,1)}
=&\qad{1,2,3}-q\qad{2,1,3}-q\qad{1,3,2}-q^3\qad{3,2,1}+q^2\qad{2,3,1}+q^2\qad{3,1,2}\\
=&(1-2\alpha q^2+2\alpha^2q^4-\alpha q^6)\det_q,
\end{align*}
then the $\Uq[3]$-cyclic span of these vectors gives $\cycmod3\alpha$
(see Example \ref{ex:hwv_for_(2,1)}).
Therefore we have
\begin{equation*}
\cycmod3\alpha \cong \begin{cases}
\Umod3{(1,1,1)} & \alpha=-1,\\
(\Umod3{(2,1)})^{\oplus 2} \oplus \Umod3{(1,1,1)} & \alpha=-1/2,\\
\Umod3{(3)}\oplus \Umod3{(2,1)}\oplus \Umod3{(1,1,1)}
& \alpha=1/(q^2\pm(q-q^3)),\\
\Umod3{(3)}\oplus (\Umod3{(2,1)})^{\oplus 2}
& \alpha=({2q^{-2}+q^2\pm\sqrt{q^4+4-4q^{-4}}})/4,\\
\Umod3{(3)}\oplus (\Umod3{(2,1)})^{\oplus 2}\oplus \Umod3{(1,1,1)}
& \text{otherwise}.
\end{cases}
\end{equation*}
In other words, we have
\begin{equation*}
\begin{split}
\mult{(3)}\alpha&=\begin{cases}
0 & \alpha=-1,-\frac12,\\
1 & \text{otherwise},
\end{cases}\\
\mult{(2,1)}\alpha&=\begin{cases}
0 & \alpha=-1,\\
1 & \alpha=1/(q^2\pm(q-q^3)),\\
2 & \text{otherwise},
\end{cases}\\
\mult{(1,1,1)}\alpha&=\begin{cases}
0 & \alpha=({2q^{-2}+q^2\pm\sqrt{q^4+4-4q^{-4}}})/4,\\
1 & \text{otherwise}.
\end{cases}
\end{split}
\end{equation*}
It also follows that
\begin{equation*}
\singular 3=\Bigl\{-1,-\frac12,
\frac1{q^2\pm(q-q^3)},\frac{{2q^{-2}+q^2\pm\sqrt{q^4+4-4q^{-4}}}}4\Bigr\}.
\end{equation*}
Notice that
$0<\mult{(2,1)}\alpha<f^{(2,1)}$ when $\alpha=1/(q^2\pm(q-q^3))$.
In particular, this shows that
the classical result \emph{cannot} be recovered
from the quantum case by letting $q \to 1$.
This is because the elementary divisors of
the $q$-content transition matrices for generic $q$
(defined in Section \ref{Classical result} below)
are different from those of
the $1$-content transition matrices;
See Example \ref{ex:F(2,1)}.
See also Example \ref{ex:n=3_case}
for the $q$-content discriminants.
\eoe\end{ex}

\section{Irreducible decomposition of $\cycmod n\alpha$
for $\alpha\in\singular n$}\label{Main}

In this section, we investigate the cases
where some irreducible factors
of the decomposition of the cyclic module $\cycmod n\alpha$
may collapse.

\subsection{$q$-content discriminants}

If $\alpha\in\singular n$,
then there exists some $\lambda\vdash n$ such that
$\mult\lambda\alpha<f^\lambda$ by definition.
To describe the sets
\begin{equation*}
\begin{split}
\subsingular n\lambda
\deq\set[1]{\alpha\in\singular n}{\mult\lambda\alpha<f^\lambda},\qquad
\singularwt n\alpha
\deq\set[1]{\lambda\in\domlattice}{\mult\lambda\alpha<f^\lambda},
\end{split}
\end{equation*}
we introduce certain polynomials called the $q$-discriminants:
Let $\lambda\in\domlattice$ and $\alpha,\,\beta\in\C$.
When $\beta\in\C\setminus\singular n$,
each vector $\hwv T$ is written as a linear combination
of the vectors $\{\hwv[\beta]T\}_{T\in\STab(\lambda)}$
\begin{equation*}
\hwv T=\sum_{S\in\STab(\lambda)}
F_q^\lambda(\alpha,\beta;{S,T})\hwv[\beta]S
\qquad(T\in\STab(\lambda)).
\end{equation*}
We introduce a $f^\lambda\times f^\lambda$ matrix
$F_{q}^\lambda(\alpha,\beta)$ by
\begin{equation*}
F_{q}^\lambda(\alpha,\beta)
=(F_q^\lambda(\alpha,\beta;{S,T}))_{S,T\in\STab(\lambda)}.
\end{equation*}
We call $F_q^\lambda(\alpha,\beta)$
the \emph{$q$-content transition matrix} of $\lambda$.
The function $C_q^\lambda(\alpha,\beta)
\deq\det F_q^\lambda(\alpha,\beta)$ of $\alpha$
is called the \emph{$q$-content discriminant} for $\lambda$
with reference point $\beta$.
If $\beta,\gamma\in\C\setminus\singular n$, then
$F_q^\lambda(\alpha,\beta)F_q^\lambda(\beta,\gamma)=F_q^\lambda(\alpha,\gamma)$
and $C_q^\lambda(\alpha,\beta)C_q^\lambda(\beta,\gamma)=C_q^\lambda(\alpha,\gamma)$.
In particular,
if $\alpha,\beta\in\C\setminus\singular n$,
then $C_q^\lambda(\alpha,\beta)=C_q^\lambda(\alpha,0)/C_q^\lambda(\beta,0)$.
In what follows,
we simply write $F_q^\lambda(\alpha)$ and $C_q^\lambda(\alpha)$
instead of $F_q^\lambda(\alpha,0)$ and $C_q^\lambda(\alpha,0)$.
By definition, we have
\begin{equation*}
\mult\lambda\alpha=\dim_\C W_{n,q}^{(\alpha)}(\lambda)=\rank F_q^\lambda(\alpha),
\end{equation*}
and hence
\begin{equation*}
\subsingular n\lambda
=\set[1]{\alpha\in\C}{C_q^\lambda(\alpha)=0}.
\end{equation*}

\begin{ex}\label{ex:F(2,1)}
Since
\begin{align*}
\hwv{\minitab(12,3)}=&(1+\alpha)(1+\alpha-2\alpha q^2)\hwv[0]{\minitab(12,3)}
+\alpha q(1+\alpha)(1-q^2)\hwv[0]{\minitab(13,2)},\\
\hwv{\minitab(13,2)}=&-\alpha q^{-1}(1+\alpha)(1-q^2)^2\hwv[0]{\minitab(12,3)}
+(1+\alpha)(1-\alpha)\hwv[0]{\minitab(13,2)},
\end{align*}
we have
\begin{equation*}
F_q^{(2,1)}(\alpha)=
\begin{pmatrix}
(1+\alpha)(1+\alpha-2\alpha q^2) & \alpha q(1+\alpha)(1-q^2)\\
-\alpha q^{-1}(1+\alpha)(1-q^2)^2 & (1+\alpha)(1-\alpha)
\end{pmatrix}
\end{equation*}
and
\begin{equation*}
\begin{split}
C_q^{(2,1)}(\alpha)&=\det F_q^{(2,1)}(\alpha)
=(1+\alpha)^2(1+(q-q^2-q^3)\alpha)(1+(-q-q^2+q^3)\alpha).
\end{split}
\end{equation*}
We note that the elementary divisors of the transition matrix
$F^{(2,1)}(\alpha)$ are given by
\begin{equation*}
\begin{cases}
(1+\alpha)(1-\alpha),\, (1+\alpha)(1-\alpha), & q^2=1,\\
(1+\alpha),\, (1+\alpha)(1+(q-q^2-q^3)\alpha)(1+(-q-q^2+q^3)\alpha), & q^2\ne1.
\end{cases}
\end{equation*}
\eoe\end{ex}

The relation between the two collections
$\{C_q^\lambda(\alpha)\}_\lambda$
and $\{\uC_q^\lambda(\alpha)\}_\lambda$ is given by the
\begin{lem}\label{lem:reduction_1}
The equalities
\begin{align*}
\rank \uF_q^\mu(\alpha)=
\sum_{\lambda\vdash n}K_{\lambda\mu} \rank F_q^\lambda(\alpha),\qquad
\uC_q^\mu(\alpha)=\prod_{\lambda\vdash n}C_q^\lambda(\alpha)^{K_{\lambda\mu}}
\end{align*}
hold where $K_{\lambda\mu}$ is the Kostka number
\pp{we refer {\upshape\cite{Mac}} for the definition}.
\end{lem}

\begin{proof}
We define
\begin{equation*}
\begin{split}
\cycmod n\alpha(\lambda,\mu)\deq
\Uq\cdot W_{n,q}^{(\alpha)}(\lambda)\cap H_\mu^{(\alpha)},\qquad
\cycmod n\alpha(T,\mu)\deq \Uq\cdot \hwv T \cap H_\mu^{(\alpha)},
\end{split}
\end{equation*}
where
\begin{equation*}
H_\mu^{(\alpha)}\deq
\sum_{\sigma\in\sym\mu\backslash\sym n}\C\cdot\qad{\vk(\mu)^\sigma}
\end{equation*}
is the subspace of $\cycmod n\alpha$
consisting of all weight vectors of weight $\mu$.
By definition, we have
\begin{equation*}
H_\mu^{(\alpha)}=\bigoplus_{\lambda\vdash n}\cycmod n\alpha(\lambda,\mu).
\end{equation*}
Since
$W_{n,q}^{(\alpha)}(\lambda)=\sum_{T\in\STab(\lambda)}\C\cdot \hwv T$,
we also have
\begin{equation*}
\cycmod n\alpha(\lambda,\mu)
=\sum_{T\in\STab(\lambda)}\cycmod n\alpha(T,\mu).
\end{equation*}
Notice that
\begin{equation*}
\cycmod n\alpha(T,\mu)\cong\begin{cases}
\Umod n\lambda(\mu) & \hwv T\ne0,\\
0 & \hwv T=0,
\end{cases}
\end{equation*}
where $\Umod n\lambda(\mu)$ is the weight space
in $\Umod n\lambda$ of weight $\mu$.
We denote by $\iota_T$ the intertwiner between
$\cycmod n\alpha(T,\mu)$
and $\Umod n\lambda(\mu)$ when $\hwv T\ne0$.
We also put
$\iota_T(x)=0\in \Umod n\lambda(\mu)$
for any $x\in \cycmod n\alpha(T,\mu)$ when $\hwv T=0$.
The map
\begin{equation*}
\cycmod n\alpha(\lambda,\mu)\supset \cycmod n\alpha(T,\mu)\ni
x \longmapsto \iota_T(x)\otimes \hwv T
\in \Umod n\lambda(\mu)\otimes W_{n,q}^{(\alpha)}(\lambda)
\end{equation*}
defines a linear isomorphism, and hence yields
\begin{align*}
\cycmod n\alpha(\lambda,\mu)
\cong \Umod n\lambda(\mu)\otimes W_{n,q}^{(\alpha)}(\lambda).
\end{align*}
Consequently, we have the decomposition
\begin{equation*}
H^{(\alpha)}_{\mu}\cong
\bigoplus_{\lambda\vdash n}\Umod n\lambda(\mu)\otimes W_{n,q}^{(\alpha)}(\lambda)
\end{equation*}
as a vector space.
Since $\dim \Umod n\lambda(\mu)=K_{\lambda\mu}$,
we have the lemma.
\end{proof}

As a corollary, we also have the
\begin{lem}\label{lem:cor_of_reduction_1}
The equalities
\begin{equation*}\label{eq:reduction_1}
\rank F_q^\lambda(\alpha)
=\sum_{\mu\vdash n}
K^{(-1)}_{\lambda\mu}\rank \uF_q^\mu(\alpha),\qquad
C_q^\lambda(\alpha)
=\prod_{\mu\vdash n}
\uC_q^\mu(\alpha)^{K^{(-1)}_{\lambda\mu}}
\end{equation*}
hold where $K^{(-1)}_{\lambda\mu}$ is the reverse Kostka number
\pp{i.e. $\sum_{\nu\vdash n}K_{\lambda\nu}K_{\nu\mu}^{(-1)}=\delta_{\lambda\mu}$}.
\qed
\end{lem}

By Lemma \ref{lem:reduction_1}, we notice that
\begin{equation*}
\uC_q^\mu(\alpha)=
C_q^\mu(\alpha)\times
\prod_{\substack{\lambda\vdash n\\ \lambda\ne\mu}}
C_q^\lambda(\alpha)^{K_{\lambda\mu}},
\end{equation*}
which readily implies that
\begin{equation*}
\singular n
=\bigcup_{\lambda\vdash n}
\set[1]{\alpha\in\C}{\uC_q^\lambda(\alpha)=0}
=\bigcup_{\lambda\vdash n}
\set[1]{\alpha\in\C}{C_q^\lambda(\alpha)=0}.
\end{equation*}
Namely, the two collections
$\{F_{q}^\lambda(\alpha)\}_{\lambda\vdash n}$ and
$\{\uF_{q}^\lambda(\alpha)\}_{\lambda\vdash n}$
of matrices have equivalent information on $\singular n$.

\subsection{Classical result --- a review}\label{Classical result}

We recall the result of the classical case \cite{MW2005}.
The set $\singular[1]n$ is explicitly given by
\begin{equation*}
\singular[1]n=\Bigl\{\pm1,\pm\frac12,\dots,\pm\frac1{n-1}\Bigr\}.
\end{equation*}
The irreducible decomposition of $\cycmod[1]n{\pm\frac1k}$
($k=1,2,\dots,n-1$) is
\begin{equation}\label{ClassicalTheorem}
\cycmod[1]n{-\frac1k}\cong\bigoplus_{%
\begin{subarray}{c}\lambda\vdash n\\ \lambda_1\le k\end{subarray}}
(\boldsymbol{E}_1^\lambda)^{\oplus f^\lambda},\qquad
\cycmod[1]n{\frac1k}\cong\bigoplus_{%
\begin{subarray}{c}\lambda\vdash n\\ \lambda'_1\le k\end{subarray}}
(\boldsymbol{E}_1^\lambda)^{\oplus f^\lambda}.
\end{equation}
In other words, the multiplicity $\mult[1]\lambda\alpha$
of the Schur module
$\boldsymbol{E}_1^\lambda$ in $\cycmod[1]n{\pm\frac1k}$ is given by
\begin{equation}\label{eq:multiplicity_in_classical_case}
\mult[1]\lambda{-\frac1k}=\begin{cases}
f^\lambda & \lambda_1\le k,\\
0 & \text{otherwise},
\end{cases}
\qquad
\mult[1]\lambda{\frac1k}=\begin{cases}
f^\lambda & \lambda'_1\le k,\\
0 & \text{otherwise}.
\end{cases}
\end{equation}
In particular, since $f^{\lambda'}=f^\lambda$,
we see that $m_1^\lambda(-\frac1k)=m_1^{\lambda'}(\frac1k)$,
where $\lambda'$ denotes the transposition of $\lambda$ as a Young diagram.
Thus, in the classical case,
each isotypic component either exists with full multiplicity or
disappears completely.
Using the (modified) content polynomial
\begin{equation*}
c^\lambda(\alpha)\deq\prod_{(i,j)\in\lambda}(1+(j-i)\alpha)
\end{equation*}
of a Young diagram $\lambda$, one finds that
\eqref{eq:multiplicity_in_classical_case}
is equivalent to
\begin{equation*}
\mult[1]\lambda\alpha=\begin{cases}
f^\lambda & c^\lambda(\alpha)\neq0,\\
0 & c^\lambda(\alpha)=0.
\end{cases}
\end{equation*}

In the quantum case,
if $\alpha\in\singular n$,
then we have
$\mult\lambda\alpha<f^\lambda$
for some $\lambda$.
In contrast to the classical case, however,
it could be that
$\mult\lambda\alpha\ne0$ as we see
in Example \ref{ex:description_for_n=3_case}.

\subsection{Symmetric and skew-symmetric cases}

The cases where $\lambda=(n)$ and $(1,\dots,1)$
are much easier
because the respective multiplicities $\mult\lambda\alpha$
of $\Umod n{\lambda}$ in the irreducible decomposition of
$\cycmod n\alpha$
are always either $0$ or $1$.
Actually, we have the following.
\begin{prop}\label{prop:extremal}
The highest weight vectors for $(n)$ and $(1,\dots,1)$
in $\cycmod n\alpha$ are
\begin{align}\label{eq:special_case}
\hwv{\minitab(12{\cdot}n)}
=\left(\sum_{\sigma\in\sym n}\alpha^{n-\cyc\sigma}\right)
\hwv[0]{\minitab(12{\cdot}n)},\qquad
\hwv{\minitab(1,2,\cdot,n)}
=\left(\sum_{\sigma\in\sym n}\alpha^{n-\cyc\sigma}(-q^2)^{\inv\sigma}\right)
\hwv[0]{\minitab(1,2,\cdot,n)}.
\end{align}
In other words, the corresponding $q$-content discriminants
\pp{and/or $q$-transition matrices} are
\begin{align*}
C_q^{(n)}(\alpha)
=F_q^{(n)}(\alpha)
=\sum_{\sigma\in\sym n}\alpha^{n-\cyc\sigma},\qquad
C_q^{(1,\dots,1)}(\alpha)
=F_q^{(1,\dots,1)}(\alpha)
=\sum_{\sigma\in\sym n}\alpha^{n-\cyc\sigma}(-q^2)^{\inv\sigma}.
\end{align*}
In particular,
\begin{equation}\label{SpecialContentPoly}
\begin{split}
\mult{(n)}\alpha=\begin{cases}
0 & \alpha=-1,-\frac12,\dots,-\frac1{n-1},\\
1 & \text{otherwise},
\end{cases}\qquad
\mult{(1,\dots,1)}\alpha=\begin{cases}
0 & \sum_{\sigma\in\sym n}\alpha^{n-\cyc\sigma}(-q^2)^{\inv\sigma}=0,\\
1 & \text{otherwise}.
\end{cases}
\end{split}
\end{equation}
\end{prop}

\begin{proof}
The first equation in \eqref{eq:special_case} is straightforward.
To prove the second one,
notice that
$\qad{1,2,\dots,n}\cdot\pi^{(\alpha)}(\yngsym(\minitab(1,2,\cdot,n)))$
is contained in the one dimensional invariant subspace $\C\cdot\det_q$,
and hence it must be a scalar multiple of $\det_q$.
The scalar is given by the coefficient of $x_{11}\dotsb x_{nn}$ in
\begin{equation*}
\begin{split}
\hwv{\minitab(1,2,\cdot,n)}
&=\qad{1,2,\dots,n}\cdot\pi^{(\alpha)}(\yngsym(\minitab(1,2,\cdot,n)))
=\sum_{w\in\sym n}(-q)^{\inv w}\qad{w(1),\dots,w(n)}\\
&=\sum_{w\in\sym n}(-q)^{\inv w}
\sum_{\sigma\in\sym n}\alpha^{n-\cyc\sigma}q^{\inv\sigma}
x_{\sigma(1)w(1)}\dotsb x_{\sigma(n)w(n)}.
\end{split}
\end{equation*}
The coefficient is equal to
\begin{equation*}
\sum_{w=\sigma\in\sym n}(-q)^{\inv w}\alpha^{n-\cyc\sigma}q^{\inv \sigma}
=\sum_{\sigma\in\sym n}\alpha^{n-\cyc\sigma}(-q^2)^{\inv \sigma}
\end{equation*}
as desired. The last statement
\eqref{SpecialContentPoly} about the multiplicity
$\mult{(n)}\alpha$ follows immediately from the fact
\begin{equation*}
\sum_{\sigma\in\sym n}\alpha^{n-\cyc\sigma}
=\prod_{k=1}^{n-1}(1+k\alpha)
\end{equation*}
(see, e.g. \cite{MW2005} or \cite{St}).
This shows the proposition.
\end{proof}

As a corollary, we have the following.
\begin{cor}
Define the set $\psingular n$ by
$\psingular n=\subsingular n{(n)}
\cup\subsingular n{(1,\dots,1)}$. Then
\begin{equation*}
\psingular n=\Bigl\{-1,-\frac12,\dots,-\frac1{n-1}\Bigr\}
\cup
\set[2]{\alpha\in\C}{\sum_{w\in\sym n}(-q^2)^{\inv w}\alpha^{n-\cyc w}=0}
\end{equation*}
and $\psingular n$ is a subset of $\singular n$.
\qed
\end{cor}
We call the singular points in $\subsingular n{(n)}$
\emph{classical} and the one in $\subsingular n{(1,\dots,1)}$
\emph{semi-classical}.
We notice that
$\psingular[1]n=\{\pm1,\pm\frac12,\dots,\pm\frac1{n-1}\}=\singular[1]n$
in the classical case.
However, it could be true
that $\singular n\supsetneq\psingular n$
in the quantum case (see, e.g. Example \ref{ex:description_for_n=3_case}).
When $\alpha\in\singular n\backslash\psingular n$, we call
it a \emph{quantum} singular point.

\subsection{Several explicit points in $\singular n$}

We use the following lemma.
\begin{lem}[{\cite[Lemma 2.1]{KW2006}}]\label{lem:partial_sum_formula_for_nu}
For any $g\in\sym n$, the equality
\begin{equation*}
\sum_{w\in\sym k}\alpha^{n-\cyc{wg}}
=\alpha^{n-\cyc{w_0g}}(1+\alpha)\dotsb(1+(k-1)\alpha)
\end{equation*}
holds.
Here $\sym k$ is regarded as a subgroup
$\sym k=\{w\in\sym n\,;\,w(x)=x,\,x>k\}$
of $\sym n$, and $w_0$ is the element in $\sym k$
\pp{depending on $g$}
such that
$\cyc{w_0g}\ge\cyc{wg}$ for any $w\in\sym k$.
\qed
\end{lem}


The following lemma plays a key role for understanding the
multiplicity $\mult\lambda{-\frac1k}$.
\begin{lem}\label{lem:vanishing_lemma}
Let $k$ be a positive integer less that $n$.
If $\alpha=-\frac1k$, then
\begin{equation}\label{k+1-vanishing}
\qad{a_1,\dots,a_m,1^{k+1},b_1,\dots,b_l}=0
\end{equation}
for any $a_1,\dots,a_m,b_1,\dots,b_l$.
Here $1^{k+1}$ denotes the $k+1$ consecutive sequence
${\overbrace{1,\dots,1}^{k+1}}$.
\end{lem}

\begin{proof}
Let $I_{n,k}$ be the set consisting of finite sequences
$(i_1,\dots,i_k)\in\{1,2,\dots,n\}^k$ such that
the entries are distinct.
We write $\vi\cap\vj=\emptyset$ if
$\vi\in I_{n,m}$ and $\vj\in I_{n,l}$ have no common entry.
For any pair $\vi=(i_1,\dots,i_m)\in I_{n,m}$ and
$\vj=(j_1,\dots,j_l)\in I_{n,l}$ such that $\vi\cap\vj=\emptyset$,
we put
\begin{equation*}
W_n(\vi,\vj)
=\set[1]{w\in\sym n}{w(x)=i_x\,(1\le x\le m),\,w(n-l+y)=j_y\,(1\le y\le l)}.
\end{equation*}
Then we have
\begin{equation}\label{eq:expansion_of_[a1b]}
\begin{split}
&\qad{a_1,\dots,a_m,1^{k+1},b_1,\dots,b_l}\\
=&\sum_{w\in\sym n}\alpha^{n-\cyc w}q^{\inv w}
x_{w(1)a_1}\dotsb x_{w(m)a_m}x_{w(m+1)1}\dotsb x_{w(m+k+1)1}
x_{w(n-l+1)b_1}\dotsb x_{w(n)b_l}\\
=&\sum_{\begin{subarray}{c}\vi=(i_1,\dots,i_m)\in I_{n,m}\\
\vj=(j_1,\dots,j_l)\in I_{n,l}\\
\vi\cap\vj=\emptyset
\end{subarray}}
\sum_{w\in W_n(\vi,\vj)}\alpha^{n-\cyc w}q^{\inv w}
x_{i_1a_1}\dotsb x_{i_ma_m}x_{w(m+1)1}\dotsb x_{w(m+k+1)1}
x_{j_1b_1}\dotsb x_{j_lb_l},
\end{split}
\end{equation}
where we suppose that $m+(k+1)+l=n$.
To prove the lemma, we show that each sum
\begin{equation}\label{eq:[a1b]-monomial}
\sum_{w\in W_n(\vi,\vj)}\alpha^{n-\cyc w}q^{\inv w}
x_{i_1a_1}\dotsb x_{i_ma_m}x_{w(m+1)1}\dotsb x_{w(m+k+1)1}
x_{j_1b_1}\dotsb x_{j_lb_l}
\end{equation}
in \eqref{eq:expansion_of_[a1b]} has a factor $(1+\alpha)\dotsb(1+k\alpha)$.
Consider the group
\begin{equation*}
\sym n(m,l)=\set[1]{w\in\sym n}{w(x)=x\ (x\le m, n-l+1\le x)}.
\end{equation*}
The group $\sym{n}(m,l)$ acts on $W_n(\vi,\vj)$ from right
transitively and faithfully.
We take the unique element $w_0\in W_n(\vi,\vj)$ such that
$w_0(x)<w_0(y)$ for $m+1\le \forall x<\forall y\le m+k+1$.
Then we have $W_n(\vi,\vj)=w_0\cdot\sym{n}(m,l)$
and $\inv{w_0w}=\inv{w_0}+\inv w$
for $w\in\sym{n}(m,l)$.
Therefore, the sum \eqref{eq:[a1b]-monomial} is rewritten as
\begin{equation*}
\begin{split}
\sum_{w\in \sym{n}(m,l)}\alpha^{n-\cyc{w_0w}}q^{\inv{w_0w}}
x_{i_1a_1}\dotsb x_{i_ma_m}x_{w_0w(m+1)1}\dotsb x_{w_0w(m+k+1)1}
x_{j_1b_1}\dotsb x_{j_lb_l}\\
=\kakko[3]{\sum_{w\in \sym{n}(m,l)}\alpha^{n-\cyc{w_0w}}}
q^{\inv{w_0}}x_{i_1a_1}\dotsb x_{i_ma_m}x_{w_0(m+1)1}\dotsb x_{w_0(m+k+1)1}
x_{j_1b_1}\dotsb x_{j_lb_l}.
\end{split}
\end{equation*}
Since $\cyc\cdot$ is a class function and
$\sym{n}(m,l)=g\sym{k+1}g^{-1}$ for some $g\in\sym n$,
we have
\begin{equation*}
\sum_{w\in \sym{n}(m,l)}\alpha^{n-\cyc{w_0w}}
=\sum_{w\in \sym{k+1}}\alpha^{n-\cyc{w_0gwg^{-1}}}
=\sum_{w\in \sym{k+1}}\alpha^{n-\cyc{w\cdot g^{-1}w_0g}},
\end{equation*}
which indeed has a factor $(1+\alpha)\dotsb(1+k\alpha)$
by Lemma \ref{lem:partial_sum_formula_for_nu}.
Thus we have the lemma.
\end{proof}

\begin{rem}
Even though there are $k+1$ identical columns,
if they are not consecutive like \eqref{k+1-vanishing},
then Lemma \ref{lem:vanishing_lemma} does not hold.
For the classical case,
such a consecutiveness condition is unnecessary
to hold the vanishment.
See \cite{KW2006}.
\end{rem}

As a corollary of the lemma above, we have the
\begin{cor}\label{cor:vanishing-cor}
If $\frac{\lambda_1}{1+\abs{\lambda}}>1-\frac1{1+k}$,
then $\mult\lambda{-\frac1k}=0$.
\end{cor}
\begin{proof}
Recall that the vectors
$\hwv T= \qad{\vj(T)}\cdot\pi^{(\alpha)}(\yngsym(T))
\in\cycmod n\alpha(\lambda)$
for $T\in\STab(\lambda)$ of size
$n$ ($\lambda\vdash n$)
form a basis of the space of highest weight vectors
$W_{n,q}^{(\alpha)}(\lambda)$ in
$\cycmod n\alpha(\lambda)$.
Here $\vj(T)=(j_1(T),\dots,j_n(T))\in I_n(\lambda)$
is defined in \eqref{eq:JT}.
Now it is easy to see that
the inequality $\frac{\lambda_1}{1+\abs{\lambda}}>1-\frac1{1+k}$
can be written as
$k(|\lambda|-\lambda_1+1)<\lambda_1$.
In other words, we have
\begin{equation*}
k\bigl(\bigl|\bigl\{p\;;\; j_p(T)\ne1\bigr\}\bigr|+1\bigr)
<\bigl|\bigl\{p\;;\; j_p(T)=1\bigr\}\bigr|.
\end{equation*}
By the pigeonhole principle,
this implies that there exists at least one $k+1$ consecutive
sequence $1^{k+1}$ in $\vj(T)$.
Hence, by Lemma \ref{lem:vanishing_lemma} the result is immediate.
\end{proof}

Practically, using Lemma \ref{lem:vanishing_lemma}, we can estimate
the multiplicity $\mult\lambda{-\frac1k}$ for each
positive integer $k$ more accurately if
the partition $\lambda$ is given explicitly.
\begin{ex}
The vectors
\begin{equation*}
\qad{1,\dots,1,\overset{j{\text{-th}}}2,1,\dots,1}\cdot
\pi^{(\alpha)}(\yngsym(\minitab(1*\cdot*,j)))
\quad(j=2,\dots,n)
\end{equation*}
generate the space of highest weight vectors
$W_{n,q}^{(\alpha)}(\lambda)$ when $\lambda=(n-1,1)$.
Since $\qad{1,\dots,1,\overset{j{\text{-th}}}2,1,\dots,1}\ne0$
holds only if $n-k\le j\le k+1$,
we have $\mult{(n-1,1)}{-\frac1k}\le\max\{2k+2-n,\,0\}$.
In particular, we have $\mult{(n-1,1)}{-\frac1k}=0$
whenever $n>2k+1$.
Moreover, suppose $k=n-1$. Then $2k+2-n=n>f^\lambda=n-1$.
Therefore, in this case, there is at least one non-zero
vector of the form
$\qad{1,\dots,1,\overset{j{\text{-th}}}2,1,\dots,1}$
which is killed by the $q$-Young symmetrizer
$\pi^{(\alpha)}(\yngsym(\minitab(1*\cdot*,j)))$.
\eoe\end{ex}

For classical singular points,
we have the following general result.
\begin{thm}\label{thm:classical sing}
It holds that
\begin{equation*}
\singularwt n{-\frac1k}
=\set[1]{\lambda\in\domlattice}{\blf{\lambda}{\e_1}>k}
\end{equation*}
for each positive integer $k$.
\end{thm}
\begin{proof}
If $\lambda_1>k$, then there exists a standard tableau
$T\in\STab(\lambda)$ such that the highest weight vector
\begin{equation*}
\hwv T
=\qad{1^{k+1},*,\dots,*}\cdot\pi^{(\alpha)}(\yngsym(T))
\end{equation*}
vanishes when $\alpha=-\frac1k$ by Lemma \ref{lem:vanishing_lemma}.
Conversely, suppose $\lambda\in \singularwt n{-\frac1k}$.
Then, since $\mult\lambda{-\frac1k}<f^\lambda$, there exists a
standard tableau $T\in\STab(\lambda)$ such that
$\hwv[-\frac1k]T=0$.
Since we have the expression
\begin{equation*}
\hwv[-\frac1k]T=\sum_{\sigma\in\sym\lambda\backslash\sym n}
Q_T^\sigma(q)D_q^{(-\frac1k)}({\vk(\lambda)^\sigma})
\end{equation*}
by \eqref{eq:vTa}
where $Q_T^\sigma(q)$ is a polynomial in $q$, the identity
$\hwv[-\frac1k]T=0$ remains true at $q=1$.
It hence follows that $\lambda\in \singularwt[1]n{-\frac1k}$,
whence $\lambda_1>k$ by \eqref{ClassicalTheorem}.
This proves the theorem.
\end{proof}

\begin{rem}
Let $\alpha=\alpha(q)\in \singular n$. Then there
is an integer $k \; (1\leq |k|\leq n-1)$ such that
$\alpha(q)\to \frac1k$ when $q\to1$. This is because the algebraic function
$\alpha(q)$ is a
root of some $q$-content discriminant and this discriminant reduces
to a content polynomial when $q\to1$.
Therefore, there is a canonical map
$\singularwt n{\alpha(q)}\to \singularwt[1]n{\frac1k}$
when $\alpha(q)\to\frac1k$, that is, $\lambda
\in \singularwt n{\alpha(q)}$
implies $\lambda\in\singularwt[1]n{\frac1k}\;(q\to 1)$.
However, the limit formula
$\lim_{q\to1}\mult\lambda{\alpha(q)}=\mult[1]\lambda{\frac1k}$
does not hold in general.
\end{rem}

By \eqref{eq:vTa}, we see that
$\mult\lambda\alpha=\dim W_{n,q}^{(\alpha)}(\lambda)
=\rank\uQ_n^\lambda(q)\uF_{q}^\lambda(\alpha)$, where the
matrix $\uQ_n^\lambda(q)$ is defined in \eqref{eq:matrix_Q}.
Since $\uF_{q}^\lambda(\alpha)$
is the identity matrix when $\alpha=0$, we have in particular
\begin{equation*}
\rank\uQ_n^\lambda(q)=\dim W_{n,q}^{(0)}(\lambda)=f^\lambda.
\end{equation*}
Therefore, we obtain the following rough estimation of the multiplicity
\begin{equation*}
f^\lambda+\rank\uF_{q}^\lambda(\alpha)-\card{\sym\lambda\backslash\sym n}
\le \mult\lambda\alpha \le
\min\{f^\lambda,\,\rank\uF_{q}^\lambda(\alpha)\}.
\end{equation*}

Since $\Det[-1]$ equals the quantum determinant $\det_q$,
we see that
$\cycmod n{-1}$ defines the one dimensional representation.
It follows that
\begin{equation*}
\rank\uF_{q}^\lambda(-1)=\begin{cases}
1 & \lambda=\e_1+\dotsb+\e_n,\\
0 & \text{otherwise}.
\end{cases}
\end{equation*}
Hence, in particular
\begin{equation*}
\singularwt n{-1}=\domlattice\setminus\{\e_1+\dotsb+\e_n\}
=\set[1]{\lambda\in\domlattice}{\blf{\lambda}{\e_1}>1}
\end{equation*}
as we stated in Theorem \ref{thm:classical sing} and
$\mult\lambda{-1}=0$ for $\lambda\in \singularwt n{-1}$.

From these observation with the facts
from Examples \ref{ex:description_for_n=2_case} and
\ref{ex:description_for_n=3_case},
we naturally reach the
\begin{conjA}\label{conj:multiplicity_survives}
A singular value
$\alpha\in\singular n$ is a quantum
\pp{i.e. $\alpha\not\in \psingular n$}
if and only if $\mult\lambda\alpha\ne0$ for any $\lambda\in\domlattice$.
\end{conjA}

\begin{rem}
The following question comes up naturally but is non-trivial:
What can one say about the relation
between the multiplicities of the roots of $q$-content discriminants
and the multiplicities of the irreducible subrepresentations
in $\cycmod n\alpha$?
\end{rem}

\subsection{$(\Uq,\Hq)$-bimodule $\UHcycmod n\alpha$}

We define the representation $\pi$ of
the algebra $\Hq$ on the space
$\cycmod n0=\bigoplus_{1\le j_1,\dots,j_n\le n}\C\cdot x_{1j_1}\dotsb x_{nj_n}$ by
\begin{equation*}
x_{1j_1}\dotsb x_{kj_k}x_{k+1j_{k+1}}\dotsb x_{nj_n}\cdot\pi(h_k)
=x_{1j_1}\dotsb x_{k+1j_k}x_{kj_{k+1}}\dotsb x_{nj_n}.
\end{equation*}
It is immediate to see that
the $\Uq$-intertwiner
$\Phi_{n,q}^{(0)}:(\C^n)^{\otimes n}\to \cycmod n0$ is also
$\Hq$-equivariant.
Hence $\pi(\Hq)$ is the commutant of $\rho(\Uq)$
in $\End(\cycmod n0)$ by Schur-Weyl duality and vice versa
\cite{J1986} (cf. \cite{Weyl}).
Therefore, as a $(\Uq,\Hq)$-bimodule,
we have the irreducible decomposition
\begin{equation}\label{eq:irred_decomp_of_V_as_bimodule}
\cycmod n0\cong
\bigoplus_{\lambda\vdash n}\Umod n\lambda\boxtimes \Hmod n\lambda,
\end{equation}
where $\Hmod n\lambda$ is the irreducible $\Hq$-module
corresponding to the partition $\lambda$.

Consider the $(\Uq,\Hq)$-cyclic module
\begin{equation*}
\UHcycmod n\alpha\deq
\rho(\Uq)\cdot{\Det}\cdot\pi(\Hq)\subset \cycmod n0.
\end{equation*}
This is the smallest $\pi(\Hq)$-invariant subspace in $\cycmod n0$
containing $\cycmod n\alpha$.
By \eqref{eq:irred_decomp_of_V_as_bimodule},
every irreducible component in $\UHcycmod n\alpha$ is of the form
$\Umod n\lambda\boxtimes \Hmod n\lambda$ with multiplicity at most one.
Hence the irreducible decomposition of
the $(\Uq,\Hq)$-bimodule $\UHcycmod n\alpha$ is given as
\begin{equation*}
\UHcycmod n\alpha
\cong \bigoplus_{\lambda\in Y_n(\alpha)}
\Umod n\lambda\boxtimes \Hmod n\lambda,
\end{equation*}
where $Y_n(\alpha)$ is a certain subset of partitions of $n$.
Therefore, if the multiplicity $\mult\lambda\alpha$ of
the irreducible representation $\Umod n\lambda$ in
$\cycmod n\alpha$ is not zero,
then the irreducible component $\Umod n\lambda\boxtimes \Hmod n\lambda$ does appear
in the irreducible decomposition of $\UHcycmod n\alpha$.
Thus Conjecture A is restated as the following form.
\begin{conjAp}
We have
$\UHcycmod n\alpha\cong
\bigoplus_{\lambda\vdash n}\Umod n\lambda\boxtimes \Hmod n\lambda$
if and only if $\alpha\in\C\setminus\psingular n$.
\end{conjAp}
We note that $\psingular[1]n=\singular[1]n$.
Hence, the conjecture can be regarded as a quantum counterpart
of the fact that
$\UHcycmod[1]n\alpha\cong
\bigoplus_{\lambda\vdash n}\Umod[1]n\lambda\boxtimes \Hmod[1]n\lambda$
if and only if $\alpha\in\C\setminus\singular[1]n$.

\subsection{Examples of $q$-content discriminants}

We here collect several examples of
$q$-content discriminants with some of their corresponding
$q$-transition matrices.

\begin{ex}\label{ex:(2,1)}
The matrix $\uF_q^{(2,1)}(\alpha)$ corresponding to
$\uC_q^{(2,1)}(\alpha)$ in the example above is given by
\begin{equation*}
\uF_q^{(2,1)}(\alpha)=
(1+\alpha)
\begin{pmatrix}
1 & \alpha q & \alpha q^2 \\
\alpha q & 1 & \alpha q(2-q^2) \\
\alpha q^2 & \alpha q(2-q^2) & 1+2\alpha-2\alpha q^2
\end{pmatrix}.
\end{equation*}
The matrix $\uF_q^{(1,1,1)}(\alpha)$ corresponding
to $\uC_q^{(1,1,1)}(\alpha)$ is also given by
\begin{equation*}
\uF_q^{(1,1,1)}(\alpha)=
\begin{pmatrix}
1 & \alpha q & \alpha q & \alpha q^3 & \alpha^2q^2 & \alpha^2q^2 \\
\alpha q & \gamma_1 & \alpha^2q^2 & \alpha q^2\gamma_3 & \alpha q^3 & \alpha q\gamma_1 \\
\alpha q & \alpha^2q^2 & \gamma_1 & \alpha q^2\gamma_3 & \alpha q\gamma_1 & \alpha q^3 \\
\alpha q^3 & \alpha q^2\gamma_3 & \alpha q^2\gamma_3 & \gamma_5 & \gamma_4 & \gamma_4 \\
\alpha^2 q^2 & \alpha q^3 & \alpha q\gamma_1 & \gamma_4 & \gamma_2 & \alpha q^2\gamma_3 \\
\alpha^2 q^2 & \alpha q\gamma_1 & \alpha q^3 & \gamma_4 & \alpha q^2\gamma_3 & \gamma_2
\end{pmatrix},
\end{equation*}
where $\gamma_1,\dots,\gamma_5$ are given by
\begin{equation*}
\begin{split}
&\gamma_1=1+\alpha-\alpha q^2,\qquad \gamma_2=1+\alpha-\alpha q^4,
\qquad \gamma_3=1+\alpha-q^2,\qquad
\gamma_4=\alpha q(2+2\alpha-2q^2-2\alpha q^2+q^4),\\
&\gamma_5=1+3\alpha+2\alpha^2-4\alpha q^2-4\alpha^2q^2
+2\alpha q^4+2\alpha^2q^4-\alpha q^6.
\end{split}
\end{equation*}
If we put
\begin{equation*}
P=\begin{pmatrix}
q^3 & q & q & -q & -q & 1 \\
q^2 & 1-q^2 & q & -1+q^2 & q & -q \\
q^2 & q & 1-q^2 & q & -1+q^2 & -q \\
1 & -q & -q & -q & -q & -q^3 \\
q & 0 & 1-q-q^2 & 0 & 1+q-q^2 & q^2 \\
q & 1-q-q^2 & 0 & 1+q-q^2 & 0 & q^2
\end{pmatrix},
\end{equation*}
then we have
\begin{equation*}
P^{-1}\uF_q^{(1,1,1)}(\alpha)P=
\diag(\alpha_{(3)},\alpha_{(2,1)}^+,\alpha_{(2,1)}^+,
\alpha_{(2,1)}^-,\alpha_{(2,1)}^-,\alpha_{(1,1,1)}),
\end{equation*}
where $\alpha_{(3)}=(1+\alpha)(1+2\alpha)$,
$\alpha_{(2,1)}^{\pm}=(1+\alpha)(1+\alpha q^2(1\pm(q-q^{-1}))$,
$\alpha_{(1,1,1)}=1-2\alpha q^2+2\alpha^2 q^4-\alpha q^6$.
\eoe\end{ex}

\begin{ex}\label{ex:n=3_case}
The $q$-content discriminants for $3$-box diagrams are given by
\begin{align*}
C_q^{(3)}(\alpha)&=(1+\alpha)(1+2\alpha),\\
C_q^{(2,1)}(\alpha)&=(1+\alpha)^2(1+(q-q^2-q^3)\alpha)(1+(-q-q^2+q^3)\alpha),\\
C_q^{(1,1,1)}(\alpha)&=1-2\alpha q^2+2\alpha^2q^4-\alpha q^6.
\end{align*}
These are obtained by the equations (Lemma \ref{lem:reduction_1})
\begin{align*}
\uC_q^{(3)}(\alpha)=C_q^{(3)}(\alpha),\qquad
\uC_q^{(2,1)}(\alpha)=C_q^{(3)}(\alpha)C_q^{(2,1)}(\alpha),\qquad
\uC_q^{(1,1,1)}(\alpha)=C_q^{(3)}(\alpha)C_q^{(2,1)}(\alpha)^2C_q^{(1,1,1)}(\alpha),
\end{align*}
and $\uC_q^\lambda(\alpha)=\det\uF_q^\lambda(\alpha)$ are calculated
from the results in Example \ref{ex:(2,1)}.
\eoe\end{ex}

\begin{ex}\label{ex:(2,2)}
Using the fact that
\begin{align*}
\yngsym(\minitab(12,34))=&
\h1+q^{-1}\h2+(q^{-1}-q)\h3+q^{-1}\h4-q^{-1}\h5-q^{-1}\h6-\h8\\
&+(q^{-2}-1)\h9+(q^{-2}-1)\h{10}
-\h{11}-q^{-2}\h{13}-q^{-2}\h{14}\\
&-q^{-1}\h{16}-q^{-1}\h{17}+(q^{-3}-q^{-1})\h{19}+q\h{20}\\
&+q\h{21}+q^{-2}\h{22}+q^2\h{23}+\h{24},\\
\yngsym(\minitab(13,24))=&
\h1-q\h2-q\h4+q^{-1}\h5+q^{-1}\h6-(q^{-1}-q)\h7\\
&-q^{-2}\h9-q^{-2}\h{10}-\h{12}+(q^{-2}-1)\h{13}+(q^{-2}-1)\h{14}\\
&-\h{15}+q^{-1}\h{16}+q^{-1}\h{17}+(q^{-1}-q^{-3})\h{19}-q^{-1}\h{20}\\
&-q^{-1}\h{21}+q^2\h{22}+q^{-2}\h{23}+\h{24},
\end{align*}
it follows that
\begin{align*}
\hwv{\minitab(12,34)}
&=(1+\alpha)(1+2\alpha+\alpha^2-3\alpha q^2-3\alpha^2q^2+2\alpha^2q^6)
\hwv[0]{\minitab(12,34)}
+\alpha q(1+\alpha)^2(1-q^2)^2\hwv[0]{\minitab(13,24)},\\
\hwv{\minitab(13,24)}
&=-\alpha q^{-1}(1-q^2)(1+\alpha)(2+\alpha-q^2-4\alpha q^2+q^4+\alpha q^4)
\hwv[0]{\minitab(12,34)}\\
&\qquad\qquad
+(1+\alpha)(1-\alpha-\alpha^2+3\alpha^2q^2-\alpha q^4-2\alpha^2q^4+\alpha q^6)
\hwv[0]{\minitab(13,24)}.
\end{align*}
Hence the $q$-transition matrix
$F_q^{(2,2)}(\alpha)$ is given by
\begin{equation*}
\begin{split}
(1+\alpha)
\begin{pmatrix}
1+2\alpha+\alpha^2-3\alpha q^2-3\alpha^2q^2+2\alpha^2q^6
& \alpha q(1+\alpha)(1-q^2)^2\\
-\alpha q^{-1}(1-q^2)(2+\alpha-q^2-4\alpha q^2+q^4+\alpha q^4)
& (1-\alpha-\alpha^2+3\alpha^2q^2-\alpha q^4-2\alpha^2q^4+\alpha q^6)
\end{pmatrix}
\end{split}
\end{equation*}
and hence the corresponding $q$-content discriminant
$C_q^{(2,2)}(\alpha)
=\det F_q^{(2,2)}(\alpha)$ is
\begin{equation*}
\begin{split}
&(1+\alpha)^2(1+(1-3q^2-q^4+q^6)\alpha-q^2(4-6q^2+q^4-q^6+q^8)\alpha^2\\
&\qquad-2q^2(1-q^2)(1-5q^2+3q^4-q^6+q^8)\alpha^3-q^2(1-q^2)^2(1-3q^2+5q^4)\alpha^4).
\end{split}
\end{equation*}
Notice that the transition matrix $F_q^{(2,2)}(\alpha)$ does become
a scalar matrix $(1-\alpha^2)I_2$ if we let $q=1$.
\eoe\end{ex}

\begin{ex}\label{ex:(3,1)}
The matrix $\uF_q^{(3,1)}(\alpha)$
corresponding to $\uC_q^{(3,1)}(\alpha)$ is given by
\begin{equation*}
(1+\alpha)
\begin{pmatrix}
(1+2\alpha) & \alpha q(1+2\alpha) & \alpha q^2(1+2\alpha) & \alpha q^3(1+2\alpha) \\
\alpha q(1+2\alpha) & 1+\alpha+\alpha q^2
& \alpha q(2+2\alpha-q^2)&  \alpha q^2(1+2\alpha)(2-q^2) \\
\alpha q^2(1+2\alpha) & \alpha q(2+2\alpha-q^2)
& 1+2\alpha+2\alpha^2 q^2(1-q^2) & \alpha q(1+2\alpha)(3-2q^2) \\
\alpha q^3(1+2\alpha) & \alpha q^2(1+2\alpha)(2-q^2)
& \alpha q(1+2\alpha)(3-2q^2) & (1+2\alpha)(1+3\alpha(1-q^2))
\end{pmatrix}.
\end{equation*}
Further, the matrix $\uF_q^{(2,2)}(\alpha)$
corresponding to $\uC_q^{(2,2)}(\alpha)$ is given by
\begin{equation*}
(1+\alpha)
\begin{pmatrix}
(1+\alpha) & \alpha q(1+\alpha) & \alpha q^2(1+\alpha)
& \alpha q^2(1+\alpha) & \alpha^2q^3(1+\alpha) & 2\alpha^2q^4 \\
\alpha q(1+\alpha) & 1+\alpha q^2 & \alpha q \gamma_1
& \alpha q \gamma_1 & \alpha^2q^2(3-q^2) & \alpha q^3\gamma_2 \\
\alpha q^2(1+\alpha) & \alpha q \gamma_1 & \gamma_4
& 2\alpha^2q^2(2-q^2) & \alpha q\gamma_3 & \alpha q^2\gamma_5 \\
\alpha q^2(1+\alpha) & \alpha q \gamma_1 & 2\alpha^2q^2(2-q^2)
& \gamma_4 & \alpha q\gamma_3 & \alpha q^2\gamma_5 \\
\alpha q^3(1+\alpha) & \alpha^2q^2(3-q^2) & \alpha q\gamma_3
& \alpha q\gamma_3 & \gamma_7 & \alpha q\gamma_6 \\
2\alpha^2 q^4 & \alpha q^3\gamma_2 & \alpha q^2\gamma_5
& \alpha q^2\gamma_5 & \alpha q\gamma_6 & \gamma_8
\end{pmatrix},
\end{equation*}
where $\gamma_1,\dots,\gamma_8$ are given by
\begin{equation*}
\begin{split}
&\gamma_1=2-q^2+\alpha q^2,\quad
\gamma_2=1+3\alpha-2\alpha q^2,\qquad
\gamma_3=1+2\alpha+q^2-\alpha q^2-q^4,\\
&\gamma_4=1+2\alpha-2\alpha q^2+\alpha q^4,\quad
\gamma_5=2+4\alpha-q^2-3\alpha q^2,\quad
\gamma_6=\alpha q(4+6\alpha-4q^2-6\alpha q^2+q^4+\alpha q^4),\\
&\gamma_7=1+\alpha+2\alpha q^2+2\alpha^2q^2-3\alpha q^4-2\alpha^2q^4+\alpha q^6,
\quad
\gamma_8=1+5\alpha+6\alpha^2-4\alpha q^2-6\alpha^2q^2-2\alpha^2q^4+2\alpha^2q^6.
\end{split}
\end{equation*}
\eoe\end{ex}

\begin{ex}
We have
\begin{equation*}
\begin{split}
C_q^{(4)}(\alpha)
&=(1+\alpha)(1+2\alpha)(1+3\alpha),\\
C_q^{(3,1)}(\alpha)
&=(1+\alpha)^3(1+2\alpha)^2(1-q^2(2-q^2)\alpha)\\
&\qquad\times(1+(1-q^4)\alpha-q^2(4-5q^2+4q^4)\alpha^2-2q^2(2-4q^2+2q^4-q^6)\alpha^3),\\
C_q^{(2,2)}(\alpha)
&=(1+\alpha)^2(1+(1-3q^2-q^4+q^6)\alpha-q^2(4-6q^2+q^4-q^6+q^8)\alpha^2\\
&\qquad-2q^2(1-q^2)(1-5q^2+3q^4-q^6+q^8)\alpha^3-q^2(1-q^2)^2(1-3q^2+5q^4)\alpha^4).
\end{split}
\end{equation*}
These are obtained by the equations (Lemma \ref{lem:reduction_1})
\begin{equation*}
\begin{split}
\uC_q^{(4)}(\alpha)=C_q^{(4)}(\alpha),\qquad
\uC_q^{(3,1)}(\alpha)=C_q^{(4)}(\alpha)C_q^{(3,1)}(\alpha),\qquad
\uC_q^{(2,2)}(\alpha)=C_q^{(4)}(\alpha)C_q^{(3,1)}(\alpha)C_q^{(2,2)}(\alpha),
\end{split}
\end{equation*}
and $\uC_q^\lambda(\alpha)=\det\uF_q^\lambda(\alpha)$ are calculated
from the results in Examples \ref{ex:(2,2)} and \ref{ex:(3,1)}.
\eoe\end{ex}

\section{Quantum $\alpha$-permanent}\label{QAP}

A theory similar to the one developed in Section \ref{QAD}
for the quantum $\alpha$-determinants can be also established
for the quantum $\alpha$-permanent defined below.
We hence close the paper by observing two examples
concerning the cyclic $\Uq$-module $\cycmodp n\alpha$
generated by a quantum $\alpha$-permanent,
and give a conjecture on a `reciprocity'
between the multiplicities of the irreducible summands
of two modules $\cycmod n\alpha$ and $\cycmodp n\alpha$.
Introducing partition functions
for the multiplicities of respective irreducible decompositions,
we close the paper by restating a certain weaker version
of the conjecture (and also Conjecture A)
in terms of the partition functions.

\subsection{Quantum $\alpha$-permanent cyclic modules
$\cycmodp n\alpha$}

We define a \emph{quantum \pp{column} $\alpha$-permanent}
by $\Per=\det_{-q^{-1}}^{(\alpha)}=\det_{q^{-1}}^{(-\alpha)}$.
Namely, we have
\begin{equation*}
\Per\deq \sum_{w\in\sym{n}}
\alpha^{n-\cyc w}(-q)^{-\inv w}x_{w(1)1}\dotsb x_{w(n)n}
\in\Aq.
\end{equation*}
We notice that $\Per[1]$ is the quantum permanent $\per_q$.
We also remark that $\Per(X)=\Per(\tX)$
as in the case of $\Det(X)$.
For convenience, we put
\begin{equation*}
\qap{j_1,\dots,j_n}\deq
\sum_{w\in\sym{n}}\alpha^{n-\cyc w}(-q)^{-\inv w}
x_{w(1),j_1}\dots x_{w(n),j_n}.
\end{equation*}
We notice that $\Per=\qap{1,2,\dots,n}$.
Let us consider the cyclic module
\begin{equation*}
\cycmodp n\alpha\deq\rho(\Uq)\cdot\Per,
\end{equation*}
and we also define $\persingular n$ and $\multp n\alpha$
similarly to $\singular n$ and $\mult n\alpha$.
By the same discussion as in the proofs of
Proposition \ref{lem:q-analogue_of_Lemma_2.2} and
Lemma \ref{lem:action_on_dets},
we have the
\begin{lem}\label{lem:action_on_pers}
We have
\begin{equation*}
\begin{split}
\rho(q^\lambda)\cdot\qap{j_1,\dots,j_n}
&= q^{\blf{\lambda}{\e_{j_1}+\dots+\e_{j_n}}}\qap{j_1,\dots,j_n},\\
\rho(e_k)\cdot\qap{j_1,\dots,j_n}
&= \sum_{l=1}^n \delta_{j_l,k+1}q_k^l(j_1,\dots,j_n)
\qap{j_1,\dots,j_{l-1},k,j_{l+1},\dots,j_n},\\
\rho(f_k)\cdot\qap{j_1,\dots,j_n}
&= \sum_{l=1}^n \delta_{j_l,k}q_k^l(j_1,\dots,j_n)
\qap{j_1,\dots,j_{l-1},k+1,j_{l+1},\dots,j_n},
\end{split}
\end{equation*}
and
\begin{equation*}
\cycmodp n\alpha=\sum_{1\le i_1,\dots,i_n\le n}
\C\cdot\qap{i_1,\dots,i_n}.
\end{equation*}
\qed
\end{lem}
We also define
\begin{equation*}
\begin{split}
\qap{j_1,\dots,j_n}\cdot\piper^{(\alpha)}(h_k)&=
\begin{cases}
\qap{j_1,\dots,j_{k+1},j_k,\dots,j_n} & j_k<j_{k+1},\\
q^{-1}\qap{j_1,\dots,j_n} & j_k=j_{k+1},\\
\qap{j_1,\dots,j_{k+1},j_k,\dots,j_n}
-(q-q^{-1})\qap{j_1,\dots,j_n} & j_k>j_{k+1}.
\end{cases}
\end{split}
\end{equation*}
for a quantum $\alpha$-permanent $\qap{j_1,\dots,j_n}$
and $h_k\in\Hq$ as in \S 3.3. We notice that, also
in the present permanent case, the result
corresponding to Proposition \ref{prop:ird_generic} holds
when $\alpha \in \C\setminus \persingular n$.

\begin{ex}\label{ex:per:description_for_n=2_case}
The highest weight vectors in $\cycmodp2\alpha$ are
\begin{equation*}
\qap{1,1}=(1-\alpha q^{-2})x_{11}x_{21},\qquad
\qap{1,2}-q\qap{2,1}=(1+\alpha)\det_q.
\end{equation*}
Hence the component $\Umod2{(2)}$ (resp. $\Umod2{(1,1)}$) does not appear
if $\alpha=-q^2$ (resp. $\alpha=1$).
Thus we have
\begin{equation*}
\multp{(2)}\alpha=\begin{cases}
0 & \alpha=q^2,\\
1 & \text{otherwise},
\end{cases}\qquad
\multp{(1,1)}\alpha=\begin{cases}
0 & \alpha=-1,\\
1 & \text{otherwise},
\end{cases}
\end{equation*}
and $\persingular 2=\{-1,q^2\}$.
\eoe\end{ex}

\begin{ex}\label{ex:per:description_for_n=3_case}
Look at the $\Uq[3]$-module $\cycmodp3\alpha$.
We can take the highest weight vectors as follows:
\begin{align*}
u^{(3)}=&
(1-2\alpha q^{-2}+2\alpha^2q^{-4}-\alpha q^{-6})x_{11}x_{21}x_{31},\\
u_1^{(2,1)}=&
(1+\alpha)(1-(q^{-1}+q^{-2}-q^{-3})\alpha)
(x_{11}x_{21}x_{32}+(1-q)x_{11}x_{22}x_{31}-qx_{12}x_{21}x_{31}),\\
u_2^{(2,1)}=&
(1+\alpha)(1-(-q^{-1}+q^{-2}+q^{-3})\alpha)
(x_{11}x_{21}x_{32}-(1+q)x_{11}x_{22}x_{31}+qx_{12}x_{21}x_{31}),\\
u^{(1,1,1)}
=&(1+\alpha)(1+2\alpha)\det_q.
\end{align*}
Therefore, we conclude that
\begin{equation*}
\cycmodp3\alpha \cong \begin{cases}
\Umod3{(3)} & \alpha=-1,\\
\Umod3{(3)}\oplus(\Umod3{(2,1)})^{\oplus 2} & \alpha=-\frac12,\\
\Umod3{(3)}\oplus \Umod3{(2,1)}\oplus \Umod3{(1,1,1)}
& \alpha=1/(q^{-2}\pm(q^{-1}-q^{-3})),\\
(\Umod3{(2,1)})^{\oplus 2}\oplus \Umod3{(1,1,1)}
& \alpha=({2q^{2}+q^{-2}\pm\sqrt{q^{-4}+4-4q^4}})/4,\\
\Umod3{(3)}\oplus (\Umod3{(2,1)})^{\oplus 2}\oplus \Umod3{(1,1,1)}
& \text{otherwise}.
\end{cases}
\end{equation*}
In other words, each multiplicity can be described as
\begin{equation*}
\begin{split}
\multp{(3)}\alpha&=\begin{cases}
0 & \alpha=(2q^{2}+q^{-2}\pm\sqrt{q^{-4}+4-4q^4})/4,\\
1 & \text{otherwise},
\end{cases}\\
\multp{(2,1)}\alpha&=\begin{cases}
0 & \alpha=-1,\\
1 & \alpha=1/(q^{-2}\pm(q^{-1}-q^{-3})),\\
2 & \text{otherwise},
\end{cases}\\
\multp{(1,1,1)}\alpha&=\begin{cases}
0 & \alpha=-1,-\frac12,\\
1 & \text{otherwise}.
\end{cases}
\end{split}
\end{equation*}
Hence, as a counterpart of the $\singular 3$ for the quantum
$\alpha$ determinant case
(Example \ref{ex:description_for_n=3_case}),
\begin{equation*}
\persingular 3=\Bigl\{-1,-\frac12,
\frac{1}{q^{-2}\pm(q^{-1}-q^{-3})},
\frac{2q^2+q^{-2}\pm\sqrt{q^{-4}+4-4q^4}}4\Bigr\}
\end{equation*}
is the corresponding singular set
$\persingular 3$ for this permanent.
\eoe\end{ex}

\subsection{Reciprocity for multiplicities
between $\singular n$ and $\persingular n$}

By the same calculation as we did
in the proof of Proposition \ref{prop:extremal},
we get
\begin{equation}\label{SpecialContentPoly_per}
\begin{split}
\multp{(n)}\alpha
=\begin{cases}
0 & \sum_{\sigma\in\sym n}\alpha^{n-\cyc\sigma}(-q^{-2})^{\inv\sigma}=0,\\
1 & \text{otherwise},
\end{cases}\quad
\multp{(1,\dots,1)}\alpha
=\begin{cases}
0 & \alpha=-1,-\frac12,\dots,-\frac1{n-1},\\
1 & \text{otherwise}.
\end{cases}
\end{split}
\end{equation}
This shows in particular that, for the permanent case,
the singular points $-\frac1k\;(1\leq k<n)$ should be also called
\emph{classical}. Moreover,
comparing \eqref{SpecialContentPoly} and \eqref{SpecialContentPoly_per},
we have the following remarkable relations
\begin{equation*}
\multp{(n)}{\alpha(q)}=\mult{(1,\dots,1)}{\alpha(q^{-1})},
\qquad
\multp{(1,\dots,1)}{\alpha(q)}=\mult{(n)}{\alpha(q^{-1})}.
\end{equation*}
Furthermore, we find the same relations in the case where $n=3$
by comparing
Examples \ref{ex:description_for_n=3_case} with
\ref{ex:per:description_for_n=3_case}
with respect to the transposition of the diagram $\lambda$.
Thus, we naturally come to expect the following `reciprocity'
(or `mirror symmetry' with respect to the reflection in the
main diagram of $\lambda$).
\begin{conjB}\label{conj:reciprocity}
{\upshape(1)}
If $\alpha(q)\in\persingular n$,
then $\alpha(q^{-1})\in\singular n$.

\noindent
{\upshape(2)}
The map
\begin{equation*}
\persingular n\ni \alpha(q) \longmapsto \alpha(q^{-1})\in\singular n
\end{equation*}
is bijective.

\noindent
{\upshape(3)}
Let $\alpha(q)\in\persingular n$.
Then the equality
\begin{equation*}
\multp\lambda{\alpha(q)}=\mult{\lambda'}{\alpha(q^{-1})}
\end{equation*}
holds for each $\lambda\in\domlattice$.
\end{conjB}

We note that the conjecture is true when $q=1$ (see Section \ref{Classical result}).
\begin{rem}
If the conjecture is true,
then it follows from Corollary \ref{cor:vanishing-cor} that
$\multp\lambda{-\frac1k}=0$ if $\frac{\lambda_1'}{1+\abs\lambda}>1-\frac1{1+k}$.
However, we notice that
there is no permanent counterpart of Lemma \ref{lem:vanishing_lemma}.
\end{rem}

Let $\partitions$ be the set of all partitions.
Define generating functions of the multiplicities
$\mult\lambda\alpha$ and $\multp\lambda\alpha$ by
\begin{align*}
\Dtheta(t,\alpha)\deq\sum_{\lambda\in\partitions}
\frac{\mult\lambda\alpha}{f^\lambda}t^{\abs\lambda}
=\sum_{n=0}^\infty\sum_{\lambda\vdash n}\frac{\mult\lambda\alpha}{f^\lambda}t^n,
\quad
\Ptheta(t,\alpha)\deq\sum_{\lambda\in\partitions}
\frac{\multp\lambda\alpha}{f^\lambda}t^{\abs\lambda}
=\sum_{n=0}^\infty\sum_{\lambda\vdash n}
\frac{\multp\lambda\alpha}{f^\lambda}t^n.
\end{align*}
We call $\Dtheta(t,\alpha)$ (resp. $\Ptheta(t,\alpha)$)
the \emph{partition function}
of the cyclic module $\cycmod n\alpha$ (resp. $\cycmodp n\alpha$).
Obviously, one has $\Dtheta[1](t,\alpha)=\Ptheta[1](t,\alpha)$.
If $\alpha\notin\bigcup_{n=1}^\infty\singular n$
(or $\alpha\notin\bigcup_{n=1}^\infty\persingular n$),
then it is readily seen that
\begin{equation*}
\Dtheta(t,\alpha)=\Ptheta(t,\alpha)
=\prod_{i=1}^\infty\frac1{1-t^i}.
\end{equation*}
When $q=1$, it is also easily calculated as
\begin{equation*}
\Dtheta[1]\kakko{t,\pm\frac1k}=
\prod_{i=1}^{k} \frac1{1-t^i}
\end{equation*}
for $k=1,2,\dots,\infty$.

In terms of these partition functions
$\Dtheta(t,\alpha)$ and $\Ptheta(t,\alpha)$,
(a slightly weaker version of) Conjectures A and B
are respectively restated as follows.
\begin{conjAB}\label{conj:partition}
{\upshape (1)}
If $\alpha$ is a quantum singular point
and $k=1,2,\dots,\infty$,
then
\begin{equation*}
\Dtheta(t,\alpha)\ne\prod_{i=1}^k\frac1{1-t^i}.
\end{equation*}

\noindent
{\upshape (2)}
The equality
\begin{equation*}
\Dtheta(t,\alpha(q^{-1}))=\Ptheta(t,\alpha(q))
\end{equation*}
holds for $\alpha(q)\in\bigcup_{n=1}^\infty\persingular n$.
\end{conjAB}
As a step for approaching the conjecture, one perhaps needs
to examine it first when the singularity is classical: 
Calculate explicitly $\Dtheta(t,-\frac1k)$ and
$\Ptheta(t,-\frac1k)$ for $k\in\Z_{>0}$.

\begin{ackn}
We thank Sho Matsumoto for his helpful comments.
\end{ackn}

\smallskip

\begin{flushleft}
Kazufumi KIMOTO\\
Department of Mathematical Science,
University of the Ryukyus.\\
Senbaru, Nishihara, Okinawa 903-0231, JAPAN.\\
\texttt{kimoto@math.u-ryukyu.ac.jp}\\

\bigskip

Masato WAKAYAMA\\
Faculty of Mathematics,
Kyushu University.\\
Hakozaki, Fukuoka 812-8518, JAPAN.\\
\texttt{wakayama@math.kyushu-u.ac.jp}
\end{flushleft}

\end{document}